\documentclass{gtart_h}  

\def\ifplaintex{\expandafter\ifx\csname documentclass\endcsname\relax}

\ifplaintex 
\else
\fi

\expandafter\ifx\csname beginpicture\endcsname\relax
\expandafter\ifx\csname documentclass\endcsname\relax
\input pictex \else
\input prepictex \input pictex \input postpictex \fi\fi

\def\gt{{\mathsurround=0pt\it $\cal G\mskip-2mu$eometry \&\ 
$\cal T\!\!$opology}}        

\def\gtp{{\mathsurround=0pt\it $\cal G\mskip-2mu$eometry \&\ 
$\cal T\!\!$opology $\cal P\!$ublications}}  

\def\lognumber#1{\def\thelognumber{#1}}
\def\volumenumber#1{\def\thevolumenumber{#1}}
\def\papernumber#1{\def\thepapernumber{#1}}
\def\volumeyear#1{\def\thevolumeyear{#1}}

\def\pagenumbers#1#2{\def\startpage{#1}\def\finishpage{#2}}
\def\published#1{\def\publishdate{#1}}
\def\proposed#1{\def\theproposer{#1}}
\def\seconded#1{\def\theseconders{#1}}
\def\received#1{\def\receiveddate{#1}}

\def\accepted#1{\def\accepteddate{#1}}
\def\asciititle#1{\def\theasciititle{#1}}

\long\def\asciiabstract#1{\long\def\theasciiabstract{#1}}
\def\asciikeywords#1{\def\theasciikeywords{#1}}

\let\\\par\let\thelognumber\relax
\let\thevolumenumber\relax\let\thepapernumber\relax
\let\thevolumeyear\relax\let\thesamplenumber\relax\let\startpage\relax
\let\finishpage\relax\let\publishdate\relax\let\receiveddate\relax
\let\reviseddate\relax\let\accepteddate\relax\let\theasciititle\relax
\let\theasciiauthors\relax
\let\theasciiabstract\relax\let\theasciikeywords\relax
\let\theasciiemail\relax\let\theshortauthors\relax\let\theshorttitle\relax

\long\def\maketitlep{   

\count0=\startpage

\gt\hfill      
\beginpicture
\setcoordinatesystem units <0.33truein, 0.33truein> point at 2.2 0.9
\setplotsymbol ({$\cal G$})
\plotsymbolspacing=9truept
\circulararc 315 degrees from 0 1 center at 0 0
\setplotsymbol ({$\cal T$})
\circulararc 315 degrees from 1 -1 center at 1 0
\endpicture
\break
{\small\ifx\thesamplenumber\relax 
Volume \else Sample
\fi\thevolumenumber\ (\thevolumeyear)
\startpage--\finishpage\nl
Published: \publishdate}
\vglue 0.5truein plus 0.4fil minus 0.1truein

{\parskip=0pt\leftskip 0pt plus 1fil\def\\{\par\smallskip}{\ifplaintex\large
\else\Large\fi\bf\thetitle}\par\medskip}   

\vglue 0pt plus 0.1fil 

{\parskip=0pt\leftskip 0pt plus 1fil\def\\{\par}{\sc\theauthors}
\par\medskip}

\vglue 0pt plus 0.1fil 

{\small\parskip=0pt\let\newline\\
{\leftskip 0pt plus 1fil\def\\{\par}{\sl\theaddress}\par}
\expandafter\ifx\theemail\relax    
\relax\else\vglue 5pt plus 0.02fil minus 2pt\def\\{\stdspace{\rm 
and}\stdspace} 
\cl{Email:\stdspace\tt\theemail}\fi
\ifx\theurl\relax                  
\relax\else\vglue 5pt plus 0.02fil minus 2pt\def\\{\stdspace{\rm 
and}\stdspace}
\cl{URL:\stdspace\tt\theurl}\fi\par}

\vglue 7pt plus 0.3fil minus 3pt

{\bf Abstract}
\vglue 5pt plus 0.1fil minus 2pt

\theabstract

\vglue 7pt plus 0.3fil minus 3pt

{\bf AMS Classification numbers}\quad Primary:\quad \theprimaryclass

Secondary:\quad \thesecondaryclass

\vglue 5pt plus 0.3fil minus 2pt

{\bf Keywords:}\quad \thekeywords

\vglue 10pt plus 0.5fil minus 5pt

{\small  Proposed: \theproposer\hfill Received: \receiveddate\nl
Seconded: \theseconders\hfill 
\ifx\reviseddate\relax                         
Accepted: \accepteddate                        
\else
Revised: \reviseddate                          
\fi}
\eject
}       

\let\maketitlepage\maketitlep
\let\maketitle\maketitlepage

\font\phead=cmsl9 scaled 950
\font=cmsl9 scaled 1050
\font\pnum=cmbx10 scaled 913
\font\lnum=cmbx10 
\font\pfoot=cmsl9 scaled 950
\font=cmsl9 scaled 1050
\ifplaintex
\headline{\vbox to 0pt{\vskip -4.5mm\line{\small\phead\ifnum
\count0=\startpage ISSN 1364-0380 (on line)
1465-3060 (printed) \hfill {\pnum\folio}\else\ifodd\count0\def\\{ }%
\ifx\theshorttitle\relax\thetitle\else\theshorttitle\fi\hfill{\pnum\folio}
\else\def\\{ and }{\pnum\folio}\hfill\ifx\theshortauthors\relax\theauthors
\else\theshortauthors\fi\fi\fi}\vss}}
\footline{\vbox to 0pt{\vglue 0mm\line{\small\pfoot\ifnum\count0=\startpage
\copyright\ \gtp\hfill\else
\gt, Volume \thevolumenumber\ (\thevolumeyear)\hfill\fi}\vss
}}
\else
\makeatletter
\def\@oddhead{{\small\lhead\ifnum\count0=\startpage ISSN 1364-0380 (on line)
1465-3060 (printed) \hfill {\lnum\number\count0}\else\ifodd\count0
\def\\{ }\ifx\theshorttitle\relax \thetitle \else\theshorttitle\fi\hfill
{\lnum\number\count0}\else\def\\{ and }{\lnum\number\count0}
\hfill\ifx\theshortauthors\relax 
\theauthors\else\theshortauthors\fi\fi\fi}}\def\@evenhead{\@oddhead}
\def\@oddfoot{\small\lfoot\ifnum\count0=\startpage\copyright\ \gtp\hfill\else
\gt, Volume \thevolumenumber\ (\thevolumeyear)\hfill\fi}
\def\@evenfoot{\@oddfoot}
\makeatother
\fi

\newwrite\gtoutfile
\long\gdef\makeheadfile{  
{\def\\{, }\def\s{ }
\immediate\openout\gtoutfile head.xxx
\immediate\write\gtoutfile{Proxy-for: \ifx\theasciiauthors\relax
\theauthors\else\theasciiauthors\fi\s<\ifx\theasciiemail\relax\theemail\else\theasciiemail\fi>}
\immediate\write\gtoutfile{\noexpand\\}
\immediate\write\gtoutfile{Authors: \ifx\theasciiauthors\relax
\theauthors\else\theasciiauthors\fi}
{\def\\{ }\immediate\write\gtoutfile{Title: \ifx\theasciititle\relax
\thetitle\else\theasciititle\fi}}
\immediate\write\gtoutfile{Subj-class: GT or SG or MG etc}
\immediate\write\gtoutfile{MSC-class: \theprimaryclass\ifx\thesecondaryclass\relax\else, \thesecondaryclass\fi}
\immediate\write\gtoutfile{Journal-ref: Geom. Topol. \thevolumenumber
(\thevolumeyear) \startpage-\finishpage}
\immediate\write\gtoutfile{Comments: Published by Geometry and Topology at}
\immediate\write\gtoutfile{\s\s http://www.maths.warwick.ac.uk/gt/GTVol\thevolumenumber/paper\thepapernumber.abs.html}
\immediate\write\gtoutfile{\noexpand\\}
\immediate\write\gtoutfile{}
\ifx\theasciiabstract\relax
\immediate\write\gtoutfile{\theabstract}\else
\immediate\write\gtoutfile{\theasciiabstract}\fi
\immediate\write\gtoutfile{}
\immediate\write\gtoutfile{\noexpand\\}
\immediate\write\gtoutfile{}
\immediate\closeout\gtoutfile}}  

\def\maketitlepage{\maketitlep\makeheadfile}
\let\maketitle\maketitlepage

\lognumber{579}
\received{15 March 2005}
\volumenumber{9}\papernumber{54}\volumeyear{2005}
\pagenumbers{2319}{2358}   
\accepted{3 December 2005}
\proposed{Benson Farb}
\seconded{Walter Neumann, Martin Bridson}
\published{21 December 2005}

\def\S{Section }

\usepackage{amsmath, amssymb,amscd}

\newtheorem {theorem}{Theorem} [section]
\newtheorem {lemma} [theorem] {Lemma}
\newtheorem {proposition} [theorem] {Proposition}

\newtheorem {corollary} [theorem] {Corollary}

\theoremstyle{definition}
\newtheorem {definition} [theorem] {Definition}

\newtheorem {remark} [theorem] {Remark}

\newtheorem {terminology} [theorem] {Terminology}

\def\R {\mathbb R}
\def\SK {\underrightarrow{\text{Ker}}\, }
\def\A {\mathcal A}
\def\P {\mathbb P}
\def\Cinf {\mathcal C_\infty}

\begin{document}

\title[Limit groups for relatively hyperbolic groups, II]{Limit groups for relatively hyperbolic groups\\II: Makanin--Razborov diagrams}
\asciititle{Limit groups for relatively hyperbolic groups, II: Makanin-Razborov diagrams}

\author{Daniel Groves}
\address{Department of Mathematics, California Institute of Technology \\
Pasadena, CA 91125, USA} 
\email{groves@caltech.edu} 

\asciikeywords{Relatively hyperbolic groups, limit groups, R-trees}


\begin{abstract}
Let $\Gamma$ be a torsion-free group which is hyperbolic relative to a
collection of free abelian subgroups.  We construct Makanin--Razborov
diagrams for $\Gamma$.  We also prove that every system of equations
over $\Gamma$ is equivalent to a finite subsystem, and a number of
structural results about $\Gamma$--limit groups.
\end{abstract}

\asciiabstract{%
Let Gamma be a torsion-free group which is hyperbolic relative to a
collection of free abelian subgroups.  We construct Makanin-Razborov
diagrams for Gamma.  We also prove that every system of equations over
Gamma is equivalent to a finite subsystem, and a number of structural
results about Gamma-limit groups.}

\maketitle

\section{Introduction}

This paper is a continuation of \cite{RelHyp}.  Throughout this paper, $\Gamma$ will denote a torsion-free group which is hyperbolic relative to a collection of free abelian subgroups.  For an arbitrary finitely generated group $G$, we wish to understand the set $\text{Hom}(G,\Gamma)$ of all homomorphisms from $G$ to $\Gamma$.  

In \cite{RelHyp} we considered a sequence of pairwise non-conjugate homomorphisms $\{ h_i \co G \to \Gamma \}$ and extracted a limiting $G$--action on a suitable asymptotic cone, and then extracted an $\R$--tree with a nontrivial $G$--action.  This $\R$--tree allows much information to be obtained.  In particular, in case $G = \Gamma$, we studied $\text{Aut}(\Gamma)$ and also proved that $\Gamma$ is Hopfian.  In this paper, we continue this study, in case $G$ is an arbitrary finitely generated group.  In particular, we construct a {\em Makanin--Razborov diagram} for $G$, which gives a parametrisation of $\text{Hom}(G,\Gamma)$ (see Section \ref{MR} below).  We build on our work from \cite{RelHyp}, which in turn builds on our previous work of \cite{CWIF} and \cite{CWIF2}. The strategy is to follow \cite[Section 1]{SelaHyp}, though there are extra technical difficulties to deal with.

To a system of equations $\Sigma$ over $\Gamma$ in finitely many variables there is naturally associated a finitely generated group $G_\Sigma$, with generators the variables in $\Sigma$, and relations the equations.  The solutions to $\Sigma$ in $\Gamma$ are in bijection with the elements of $\text{Hom}(G_\Sigma,\Gamma)$.  Thus, Makanin--Razborov diagrams give a description of the set of solutions to a given system of equations over $\Gamma$.  For free groups, building on the work of Makanin and Razborov, Makanin--Razborov diagrams were constructed by Kharlampovich and Miasnikov \cite{KM1, KM2}, and also by Sela \cite{Sela1}.  For torsion-free hyperbolic groups, Makanin--Razborov diagrams were constructed by Sela \cite{SelaHyp}, and it is Sela's approach that we follow here.  Alibegovi\'c \cite{Alibegovic2} constructed Makanin--Razborov diagrams for limit groups.

Limit groups are hyperbolic relative to their maximal non-cyclic abelian subgroups (see \cite{Dahmani} and \cite{Alibegovic}).  Limits groups are also torsion-free.  Therefore, the main result of this paper (the construction of Makanin--Razborov diagrams) generalises the main result of \cite{Alibegovic2}.  Alibegovi\'c has another approach to the construction of Makanin--Razborov diagrams for these relatively hyperbolic groups (see \cite[Remark 3.7]{Alibegovic2}).

The main results of this paper are the following five:

\medskip

{\bf \noindent Theorem \ref{FullyResiduallyGamma}}\qua 
{\sl Suppose that $\Gamma$ is a torsion-free relatively hyperbolic group with abelian parabolics, and that $G$ is a finitely generated group.  Then $G$ is a $\Gamma$-limit group if and only if $G$ is fully residually $\Gamma$.
}

\medskip

{\bf \noindent Proposition \ref{CountablyMany}}\qua 
{\em Suppose that $\Gamma$ is a torsion-free relatively hyperbolic group with abelian parabolics.  Then there are only countably many $\Gamma$--limit groups.
}

\medskip

{\bf\noindent Theorem \ref{FiniteSubsystem}}
{\sl Let $\Gamma$ be a torsion-free relatively hyperbolic group with abelian parabolics.  Then every system of equations in finitely many
variables over $\Gamma$ (without coefficients) is equivalent to a finite subsystem.
}

\medskip

In Section \ref{ShortenSection} we define an equivalence relation
`$\sim$' on the set of homomorphisms $\text{Hom}(G,\Gamma)$, where
$G$ is
an arbitrary finitely generated group and $\Gamma$ is a torsion-free 
relatively hyperbolic group with abelian parabolics.  This
equivalence relation uses conjugation, elements of $\text{Mod}(G)
\le \text{Aut}(G)$ and `bending' moves (see Section 
\ref{ShortenSection}). 

\medskip

\begin{theorem} \label{FactorSet}
Let $\Gamma$ be a torsion-free relatively hyperbolic group with abelian parabolics, and let $G$ be a finitely generated freely indecomposable group.  There is a finite collection $\{ \eta_i \co G \to L_i \}_{i=1}^n$ of proper quotients of $G$ such that, for any homomorphism $h \co G \to \Gamma$ which is not equivalent to an injective homomorphism
there exists $h' \co G \to \Gamma$ with 
$h \sim h'$, $i \in \{ 1 , \ldots , n \}$ and $h_i \co L_i \to \Gamma$ so that $h' = h_i \circ \eta_i$.
\end{theorem}

The quotient groups $L_i$ that appear in Theorem \ref{FactorSet} are 
{\em $\Gamma$--limit groups}.  Theorem \ref{FactorSet} allows us to reduce the description of $\text{Hom}(G,\Gamma)$ to a description of $\{ \text{Hom}(L_i,\Gamma) \}_{i=1}^n$.  We then apply Theorem \ref{FactorSet} to each of the $L_i$ in turn, and so on with successive proper quotients.  That this process terminates follows from the following

\medskip

{\bf\noindent Theorem \ref{SeriesTerminates}}\qua
{\sl Let $\Gamma$ be a torsion-free group which is hyperbolic relative to free abelian subgroups.  Every decreasing sequence of $\Gamma$--limit groups:
\[	R_1 > R_2 > R_3 > \ldots	,	\]
which are all quotients of a finitely generated group $G$, 
terminates after finitely many steps.
}

\medskip

Theorems \ref{FactorSet} and \ref{SeriesTerminates} allow us to construct {\em Makanin--Razborov diagrams} over $\Gamma$, which is one of the main purposes of this paper.  A Makanin--Razborov diagram for $G$ is a finite tree which encodes all of the information about $\text{Hom}(G,\Gamma)$ obtained from the above process.  Thus a Makanin--Razborov diagram gives a parametrisation of $\text{Hom}(G,\Gamma)$.  This is described in further detail in Section \ref{MR}, the main result of which is the existence of Makanin--Razborov diagrams over $\Gamma$.

Sela \cite[I.8]{SelaProblems} asked whether Theorem \ref{FiniteSubsystem} holds and whether Makanin--Razborov diagrams exist in the context of CAT$(0)$ groups with isolated flats.  We believe that relatively hyperbolic groups with abelian parabolics are a natural setting for these questions.  However, it is worth noting that for technical reasons we have restricted to abelian, rather than virtually abelian, parabolics, so we have not entirely answered Sela's questions.

An outline of this paper is as follows. In Section \ref{Prelims} we recall the definition of relatively hyperbolic groups, and recall the construction of limiting $\R$--trees from \cite{CWIF} and \cite{RelHyp}, as well as other useful results. In Section \ref{ShortenSection} we improve upon our version of Sela's shortening argument from \cite{CWIF2} and \cite{RelHyp} to deal with arbitrary sequences of homomorphisms $\{ h_n \co G \to \Gamma \}$ where $G$ is an arbitrary finitely generated group.  In Section \ref{ShortenQSection} we recall Sela's construction of {\em shortening quotients} from \cite{Sela1}, and adapt this construction to our setting. In Section \ref{MainSection} we prove Theorem \ref{SeriesTerminates}, one of the main technical results of this paper.  We also prove Theorems \ref{CountablyMany}, \ref{FiniteSubsystem}, and a number of structural results about $\Gamma$--limit groups.  Finally in Section \ref{MR} we construct Makanin--Razborov diagrams over $\Gamma$.

\smallskip

{\bf Acknowledgment}\qua
I would like to thank Zlil Sela for providing me with the proof of \cite[Proposition 1.21]{SelaHyp}, which is repeated in the proof of Proposition \ref{FiniteMax} in this paper.  I would also like to thank
the referee for correcting a number of mistakes in earlier versions
of this paper, in particular the use of the bending moves in shortening quotients, and for his/her careful reading(s) and numerous comments, 
which have substantially improved the exposition of the results in this paper.

This work was supported in part by NSF Grant
DMS-0504251.

\section{Preliminaries} \label{Prelims}

\subsection{Relatively hyperbolic groups}

{\em Relatively hyperbolic} groups were first defined by Gromov in his seminal paper on hyperbolic groups \cite[Subsection 8.6, p.256]{Gromov}.  Another definition was given by Farb \cite[Section 3]{Farb}, and further definitions given by Bowditch \cite[Definitions 1 and 2, page 1]{Bowditch}.  These definitions are all equivalent (see \cite[Theorem 7.10, page 44]{Bowditch} and \cite[Theorem 6.1, page 682]{Dahmani2}).  Recently there has been a large amount of interest in these groups (see \cite{Alibegovic, Dahmani2, Dahmani3, DS, DS--RD, Osin, Yaman}, among others).  The definition we give here is due to Bowditch \cite[Definition 2, page 1]{Bowditch}.

\begin{definition}
A group $\Gamma$ with a family $\mathcal G$ of finitely generated subgroups is called {\em hyperbolic relative to $\mathcal G$} if $\Gamma$ acts on a $\delta$--hyperbolic graph $\mathcal K$ with finite quotient and finite edge stabilisers, where the stabilisers of infinite valence vertices are the elements of $\mathcal G$, so that $\mathcal K$ has only finitely many orbits of simple loops of length $n$ for each positive integer $n$.

The groups in $\mathcal G$ are called {\em parabolic subgroups of $\Gamma$}.
\end{definition}

In this paper we will be exclusively interested in relatively hyperbolic groups which are torsion-free and have abelian parabolic subgroups.

We record the following simple lemma for later use.
\begin{lemma} \label{Abfg}
Suppose that $\Gamma$ is a torsion free relatively hyperbolic group with abelian parabolics.  Then abelian subgroups of $\Gamma$ are finitely generated.
\end{lemma}
\begin{proof}
Suppose that $g \in \Gamma$.  If $g$ is a hyperbolic element, then by a result of Osin (see \cite[Theorem 1.14, page 10]{Osin} and the comment thereafter), the centraliser of $g$ is cyclic.

Therefore, any noncyclic abelian subgroup is contained entirely within a single parabolic subgroup.  These are assumed to be abelian (and are finitely generated by the definition of `relatively hyperbolic').
\end{proof}

\subsection{The limiting $\R$--tree} \label{LimitingSection}

In this subsection we recall a construction from \cite{RelHyp} (see also \cite{CWIF} for more details).  Suppose that $\Gamma$ is a torsion-free relatively hyperbolic group with abelian parabolic subgroups.  In \cite[Section 4]{RelHyp}, we constructed a space $X$ on which $\Gamma$ acts properly and cocompactly by isometries.  For each parabolic subgroup $P$ (of rank $n$, say) there is in $X$ an isometrically embedded copy of $\R^n$, with the Euclidean metric, so that the action of $P$ leaves this Euclidean space invariant and this $P$--action is proper and cocompact with quotient the $n$--torus.

\begin{terminology}
We say that two homomorphisms $h_1, h_2 \co G \to \Gamma$ are {\em conjugate} if there exists $\gamma \in \Gamma$ so that $h_1 = \tau_\gamma \circ h_2$, where $\tau_\gamma$ is the inner automorphism of $\Gamma$ induced by $\gamma$.  Otherwise, $h_1$ and $h_2$ are {\em non-conjugate}.
\end{terminology}

Suppose now that $G$ is a finitely generated group, and that $\{ h_n \co G \to \Gamma \}$ is a sequence of pairwise non-conjugate homomorphisms.  By considering the induced actions of $G$ on $X$, and passing to a limit, we extract an isometric action of $G$ on an asymptotic cone $X_\omega$ of $X$.  This action has no global fixed point.  There is a separable $G$--invariant subset $\mathcal C_\infty \subseteq X_\omega$, and by passing to a subsequence $\{ f_i \}$ of $\{ h_i \}$ we may assume that the (appropriately scaled) actions of $G$ on $X$ converge in the $G$--equivariant Gromov--Hausdorff topology to the $G$--action on $\mathcal C_\infty$.  

The space $\mathcal C_\infty$ is a {\em tree-graded} metric space, in the terminology of Dru\c{t}u and Sapir \cite[Definition 1.10, page 7]{DS}.  Informally, this means that there is a collection of `pieces' (in this case finite dimensional Euclidean spaces), and otherwise the space is `tree-like' (see \cite{DS} for the precise definition and many properties of tree-graded metric spaces).  By carefully choosing lines in the `pieces', and projecting, an $\R$--tree $T$ is extracted from $\mathcal C_\infty$.  This tree $T$ comes equipped with an isometric $G$--action with no global fixed points and the kernel of the $G$--action on $T$ is the same as the kernel of the $G$--action on $\mathcal C_\infty$.  For more details on this entire construction, see \cite[Sections 4,5,6]{RelHyp} and \cite{CWIF}.

\begin{definition} \cite[Definition 1.5, page 3]{BFSela}\label{Stable}\qua
Suppose that $\{ h_n \co G \to \Gamma \}$ is a sequence of homomorphisms.  The {\em stable kernel} of $\{ h_n \}$, denoted $\SK (h_n)$, is the set of all $g \in G$ so that $g \in \text{ker}(h_n)$ for all but finitely many $n$.

A sequence of homomorphisms $\{ h_i \co G \to \Gamma \}$ is {\em stable} if for all $g \in G$ either (i) $g \in \text{ker}(h_i)$ for all but finitely many $i$; or (ii) $g \not\in \text{ker}(h_i)$ for all but finitely many $i$.
\end{definition}

The following theorem recalls some of the properties of the $G$--action on the $\R$--tree $T$.

\begin{theorem} [{Compare \cite[Theorem 4.4]{CWIF} and \cite[Theorem 6.4]{RelHyp}}] \label{Texists}
Suppose that $\Gamma$ is a torsion-free group that is hyperbolic relative to a collection of free abelian subgroups and that $G$ is a finitely generated group.  Let $\{ h_n \co G \to \Gamma \}$ be a sequence of pairwise non-conjugate homomorphisms.  There is a subsequence $\{ f_i \}$ of $\{ h_i \}$ and an action of $G$ on an $\R$--tree $T$ so that if $K$ is the kernel of the $G$--action on $T$ and $L := G/K$ then:
\begin{enumerate}
\item The stabiliser in $L$ of any non-degenerate segment in $T$ is free abelian.
\item If $T$ is isometric to a real line then $L$ is free abelian, and for all but finitely many $n$ the group  $h_n(G)$ is free abelian.
\item If $g \in G$ stabilises a tripod in $T$ then $g \in \SK (f_i) \subseteq K$.
\item Let $[y_1, y_2] \subset [y_3, y_4]$ be non-degenerate segments in $T$, and assume that $\text{Stab}_L([y_3,y_4])$ is nontrivial.  Then
\[	\text{Stab}_L([y_1,y_2]) = \text{Stab}_L([y_3,y_4]).	\]
In particular, the action of $L$ on $T$ is stable.
\item $L$ is torsion-free.
\end{enumerate}
\end{theorem}

Thus $T$ is isometric to a line if and only if $L$ is abelian.  If $L$ is not abelian then $K = \SK (f_i)$.

We now recall the definition of {\em $\Gamma$--limit groups}.  There are many ways of defining $\Gamma$--limit groups.  We choose a geometric definition using the above construction.

\begin{definition}(Compare \cite[Definition 1.11]{SelaHyp}, \cite[Definition 1.2]{RelHyp})\qua
A {\em strict $\Gamma$--limit group} is a quotient $G/K$ where $G$ is a finitely generated group, and $K$ is the kernel of the $G$--action on $T$, where $T$ is the $\R$--tree arising from a sequence of non-conjugate homomorphisms $\{ h_n \co G \to \Gamma \}$ as described above.

A $\Gamma$--limit group is a group which is either a strict $\Gamma$--limit group or a finitely generated subgroup of $\Gamma$.
\end{definition}

\begin{definition}
We say that a sequence of homomorphisms $\{ h_n \co G \to \Gamma 
\}$ {\em converges to a $\Gamma$--limit group $L$} if the subsequence
in Theorem \ref{Texists} can be taken to be $\{ h_n \}$ itself and $L$
is the resulting $\Gamma$--limit group.
\end{definition}

\begin{remark}
There are finitely generated subgroups of torsion-free hyperbolic groups which are not finitely presented (see, for example, the construction from \cite[Corollary (b), page 46]{Rips}).  Therefore, when $\Gamma$ is a torsion-free relatively hyperbolic group with free abelian parabolic subgroups, a $\Gamma$--limit group need not be finitely presented.  This presents substantial complications (many of which are already dealt with by Sela in \cite[Section 1]{SelaHyp}), some of which are solved by the application of Theorem \ref{GeneralFiniteRes} below.
\end{remark}

\subsection{Acylindrical accessibility and JSJ decompositions}

In \cite{SelaAcyl}, Sela studied {\em acylindrical} graph of groups decompositions, and proved an accessibility theorem for $k$--acylindrical splittings \cite[Theorem 4.1]{SelaAcyl}.  Unlike other accessibility results such as \cite[Theorem 5.1, page 456]{Dunwoody} and \cite[Main theorem, page 451]{BF}, Sela's result holds for finitely generated groups, rather than just for finitely presented groups.\footnote{Since $\Gamma$--limit 
groups needn't be finitely presented, it is important that we can
apply acylindrical accessibility.}

We follow \cite[Section 2]{Sela1} to construct the abelian JSJ decomposition
of a freely indecomposable strict $\Gamma$--limit group, $L$.

Recall the following:
\begin{lemma} \label{AbMal}
{\rm\cite[Lemma 6.7]{RelHyp}}\qua
Let $\Gamma$ be a torsion-free group which is hyperbolic relative
to abelian subgroups, and let $L$ be a strict $\Gamma$--limit group.
Abelian subgroups of $L$ are malnormal.  Each abelian
subgroup of $L$ is contained in a unique maximal abelian 
subgroup.
\end{lemma}

Given Lemma \ref{AbMal}, the proof of the following is identical
to that of \cite[Lemma 2.1]{Sela1}.

\begin{lemma} \label{MakeAbEll}
Let $\Gamma$ be a torsion-free group which is hyperbolic relative
to abelian subgroups, and let $L$ be a strict $\Gamma$--limit group.
Let $M$ be a maximal abelian subgroup
of $L$, and let $A$ be any abelian subgroup of $L$. 
\begin{enumerate}
\item If $L = U \ast_A V$ then $M$ can be conjugated into
either $U$ or $V$.
\item If $L = U \ast_A$ then either (i) $M$ can be conjugated into
$U$ or (ii) $M$ can be conjugated to $M'$ and $L = U \ast_A M'$.
\end{enumerate}
\end{lemma}

Given Lemma \ref{MakeAbEll}, if we have an abelian splitting
$L = U \ast_A$ where $A$ is a subgroup of a non-elliptic maximal
abelian subgroup $M$, we can convert it into the amalgamated
free product $L = U \ast_A M$.
Thus we will concentrate only on splittings where edge groups are abelian and every noncyclic abelian subgroup is elliptic.
Just as in \cite[Lemma 2.3]{Sela1}, any such splitting can be modified
by sliding operations and modifying boundary monomorphisms by
conjugation to be $2$--acylindrical.  Therefore, we can apply acylindrical
accessibility.

We now follow the analysis from \cite[Section 2]{Sela1}.  First we construct
the canonical quadratic decomposition of $L$ (see \cite[Theorem 5.6, pages 98--99]{RipsSelaJSJ}.  This can then
be refined to construct the cyclic JSJ decomposition, and
further refined to construct the abelian JSJ decomposition, which
encodes all $2$--acylindrical abelian splittings of $L$ in which noncyclic
abelian subgroups are elliptic.  The proof of
\cite[Theorem 1.7]{SelaGAFA} applies in this context and implies that 
the abelian JSJ  decomposition of $L$ is unique up to sliding moves, 
conjugation and modifying boundary monomorphisms by conjugation.  
See \cite[Theorem 1.10]{SelaHyp} for a precise statement (which holds
in our context without change), and 
\cite{RipsSelaJSJ} for definitions of terms left undefined here.

There are three types of vertex groups in the abelian JSJ decomposition
of $L$.  The first are {\em surface vertex groups}, which correspond to
maximal quadratically hanging subgroups.  The second are {\em
abelian vertex groups} which are abelian. The third type are
{\em rigid vertex groups}.  The key feature of rigid vertex groups is
that they are elliptic in any abelian splitting of $L$.

By Theorem \ref{Texists} any strict $\Gamma$--limit group $L$ admits a
nontrivial stable action on any $\R$--tree with abelian segment 
stabilisers and trivial tripod stabilisers.  If $L$
is nonabelian then this $\R$--tree is not a line and, by 
\cite[Theorem 3.1]{SelaAcyl}, $L$ splits over a group of 
the form $E$--by-cyclic,
where $E$ stabilises a segment.  Since abelian subgroups of $L$
are malnormal (Lemma \ref{AbMal}) and segment stabilisers are
abelian, any group of the form $E$--by-cyclic is in fact abelian.  Therefore, any non-abelian strict $\Gamma$--limit group admits a
nontrivial abelian splitting.  This now implies the following result.

\begin{proposition} \label{JSJnontriv}
Suppose that $L$ is a non-abelian, freely indecomposable strict
$\Gamma$--limit group, and that $L$ is not a surface group.  The
abelian JSJ decomposition of $L$ is nontrivial.
\end{proposition}

\section{The shortening argument} \label{ShortenSection}

In \cite[Sections 3--6]{CWIF2} and \cite[Section 7]{RelHyp} we described a version of Sela's shortening argument which worked for sequences of {\em surjective} homomorphisms to $\Gamma$, and described in \cite[Section 3]{CWIF2} why this notion is insufficient for all sequences of homomorphisms.

In this section we present another version of the shortening argument, which works for all sequences of homomorphisms $\{ h_n \co G \to \Gamma \}$, for {\em any} finitely generated group $G$.  This version was stated but not proved in \cite[Theorem 7.10]{RelHyp}, and we give the proof here.

There are two equivalent approaches to this version of the shortening argument.  The first is to find a group $\widehat{G}$ which contains $G$ and shorten using elements of $\text{Mod}(\widehat{G})$, rather than just elements of $\text{Mod}(G)$ (this approach was used in the proof of \cite[Theorem 7.9]{CWIF2}) .  The second approach is to use the `bending' moves of Alibegovi\'c \cite{Alibegovic2}.  We use the second approach here.

\begin{definition}
Let $G$ be a finitely generated group.  A {\em Dehn twist} is an automorphism of one of the following two types:
\begin{enumerate}
\item Suppose that $G = A \ast_C B$ and that $c$ is contained in the centre of $C$.  Then define $\phi \in \text{Aut}(G)$ by $\phi(a) = a$ for $a \in A$ and $\phi(b) = cbc^{-1}$ for $b \in B$.
\item Suppose that $G = A \ast_C$, that $c$ is in the centre of $C$, and that $t$ is the stable letter of this HNN extension.  Then define $\phi \in \text{Aut}(G)$ by $\phi(a) = a$ for $a \in A$ and $\phi(t) = tc$.
\end{enumerate}
\end{definition}

\begin{definition} [Generalised Dehn twists]
Suppose $G$ has a graph of groups decomposition with abelian edge groups, and $A$ is an abelian vertex group in this decomposition.  Let $A_1 \le A$ be the subgroup generated by all edge groups connecting $A$ to other vertex groups in the decomposition.  Any automorphism of $A$ that fixes $A_1$ element-wise can be naturally extended to an automorphism of the ambient group $G$.  Such an automorphism is called a {\em generalised Dehn twist} of $G$.
\end{definition}

\begin{definition} \label{Mod}
Let $G$ be a finitely generated group.  We define ${\rm Mod}(G)$ to be the subgroup of ${\rm Aut}(G)$ generated by:
\begin{enumerate}
\item Inner automorphisms,
\item Dehn twists arising from splittings of $G$ with abelian edge groups,
\item generalised Dehn twists arising from graph of groups decompositions of $G$ with abelian edge groups.
\end{enumerate}
\end{definition}
Similar definitions are made in \cite[Section 5]{Sela1} and \cite[Section 1]{BFSela}.

We will try to shorten homomorphisms by pre-composing with elements of $\text{Mod}(G)$.  However, as seen in \cite[\S 3]{CWIF2}, this is not sufficient to get the most general result.  Thus, we also define a further kind of move (very similar to Alibegovi\'c's {\em bending} move, \cite[\S 2]{Alibegovic2}).

\begin{definition} \label{Bending}
Suppose that $\Gamma$ is a torsion-free group which is hyperbolic relative to free abelian subgroups, that $G$ is a finitely generated group, and that $h \co G \to \Gamma$ is a homomorphism.  We define two kinds of `bending' moves as follows:
\begin{enumerate}
\item[(B1)] Let $\Lambda$ be a graph of groups decomposition of $G$, and let $A$ be an abelian vertex group in $\Lambda$.  Suppose that $h(A) \le P$, where $P$ is a parabolic subgroup of $\Gamma$.  A move of type (B1) replaces $h$ by a homomorphism $h' \co G \to \Gamma$ which is such that:
\begin{enumerate}
\item $h'(A) \cong h(A)$;
\item $h'(A) \le P$;
\item $h'$ agrees with $h$ on all edge groups of $\Lambda$ adjacent to $A$, and all vertex groups other than $A$.
\end{enumerate}
\item[(B2)]  Let $\Lambda$ be a graph of groups decomposition of $G$, and let $A$ be an abelian edge group of $\Lambda$, associated to an edge $e$.  Suppose that $h(A) \le P$, for some parabolic subgroup $P$ of $\Gamma$.  A move of type (B2) replaces $h$ by a homomorphism which either
\begin{enumerate}
\item conjugates the image under $h$ of a component of $\pi_1(\Lambda \smallsetminus e)$ by an element of $P$, in case $e$ is a separating edge, or
\item multiplies the image under $h$ of the stable letter associated to $e$ by an element of $P$, in case $e$ is non-separating,
\end{enumerate}
and otherwise agrees with $h$.
\end{enumerate}
\end{definition}

\begin{definition}(Compare \cite[Definition 4.2]{BFSela}, \cite[Definition 2.11]{Alibegovic2})\label{ShortHomo}\qua
We define the relation `$\sim$' on the set of homomorphisms $h \co G \to \Gamma$ to be the equivalence relation generated by setting $h_1 \sim h_2$ if $h_2$ is obtained from $h_1$ by:
\begin{enumerate}
\item pre-composing with an element of $\text{Mod}(G)$,
\item post-composing with an inner automorphism of $\Gamma$ or
\item a bending move of type (B1) or (B2).
\end{enumerate}
\end{definition}

\begin{definition}
Let $\A$ be an arbitrary finite generating set for $G$, and let $X$ be the space upon which $\Gamma$ acts properly, cocompactly and isometrically (defined in \cite{RelHyp} and briefly described in Section 
\ref{Prelims}), with basepoint $x$.  For a homomorphism $h \co G \to \Gamma$ define $\| h \|$ by
\[	\| h \| :=  \max_{g \in \A} d_X(x , h(g) . x)	.	\]
A homomorphism $h \co G \to \Gamma$ is {\em short} if for any $h'$ such that $h \sim h'$ we have $\| h \| \le \| h' \|$.
\end{definition}

The following is the main result of this section.

\begin{theorem} \label{Short}
Suppose that $\Gamma$ is a torsion-free group which is hyperbolic relative to a collection of free abelian subgroups.  Let $G$ be a freely indecomposable finitely generated group and $\{ h_n \co G \to \Gamma \}$ be a sequence of non-conjugate homomorphisms which converges to a faithful action of $G$ on the tree-graded metric space $\Cinf$ as in Subsection \ref{LimitingSection} above.  Then, for all but finitely many $n$, the homomorphism $h_n$ is not short.
\end{theorem}
\begin{proof}
Suppose that the sequence $\{ h_n \co G \to \Gamma \}$ converges into a faithful $G$--action on $\Cinf$.  As described in Subsection \ref{LimitingSection} (see \cite{CWIF} and \cite{RelHyp} for more details), a faithful $G$--action on an $\R$--tree $T$ may be extracted from the $G$--action on $\Cinf$.  We now analyse the action of $G$ on $T$ more closely.

The group $G$ is freely indecomposable and the stabiliser in $G$ of any tripod in $T$ is trivial, so we can apply the decomposition theorem of Sela -- \cite[Theorem 3.1]{SelaAcyl} -- and decompose $T$ into subtrees of three types: axial, IET, and discrete (because $G$ is freely indecomposable and tripod stabilisers are trivial, there are no Levitt-type components in $T$).  This decomposition of $T$ induces (via the Rips machine; see \cite{BF-RipsMachine}) a graph of groups decomposition of $G$, which will allow us to shorten $h_n$ for sufficiently large $n$.  See \cite{RipsSelaGAFA}, \cite[Sections 2,3]{SelaAcyl} or \cite[Section 4]{CWIF2} for more information.

Note that there are two sources for segments in $T$.  There are segments in $\Cinf$ meeting each flat in at most a point, and there are flats in $\Cinf$ which are projected to lines in $T$.  We treat these as two separate cases.  However, we can make the following simplifications (let $\P$ be the collection of lines in $T$ which are projections of flats in $\Cinf$):
\begin{enumerate}
\item Suppose that $Y$ is an IET subtree of $T$ and that $p_E \in \P$ is a line in $T$.  Then the intersection $Y \cap p_E$ contains at most a point (\cite[Proposition 4.3]{CWIF2}).
\item  Suppose that a line $l \subset T$ is an axial subtree of $T$ and that the line $p_E$ is in $\P$.  If $l \cap p_E$ contains more than a point then $l = p_E$ (\cite[Proposition 4.5]{CWIF2}).
\item If an edge $e$ in the discrete part of $T$ has an intersection of positive length with $p_E \in \P$ then $e \subset p_E$ (\cite[Lemma 4.7]{CWIF2}).
\end{enumerate}

Fix a finite generating set $\A$ for $G$. Let $y$ be the basepoint in $T$, and consider the paths $[y, u.y]$ for $u \in \A$.  If there is any IET component of $T$ which intersects any segment $[y,u.y]$ nontrivially then we can apply \cite[Theorem 5.1]{RipsSelaGAFA} and \cite[Corollary 4.4]{CWIF2} to shorten these intersections whilst leaving the remaining segments unchanged (to see that we can have $G$ finitely generated rather than finitely presented, see \cite[Remark 4.8]{CWIF2}).

Suppose that some segment $[y,u.y]$ has an intersection of positive length with some axial component $l \subset T$ so that $l$ is not contained in any $p_E \in \P$.  Then \cite[Theorem 5.1]{CWIF2} can be used to shorten those segments $[y,u_i.y]$ intersecting the orbit of $p_E$ nontrivially, and leaving other segments unchanged.

Suppose that $[y,u.y]$ intersects some line $p_E$ nontrivially, and that $p_E$ is an axial component of $T$.  The only place where the proof of \cite[Theorem 5.2]{CWIF2} breaks down is that the images $h_i(G)$ may not intersect parabolics in its image in denser and denser subsets (when measured with the scaled metric).  However, this is exactly what the bending move (B1) is designed to deal with.  

We have the following analogue of \cite[Proposition 5.4]{CWIF2}:  Let $E$ denote the flat in $\Cinf$ which projects to $p_E$.  The subgroup $\text{Stab}_G(E)$ is an abelian subgroup of $G$.  There is a sequence of flats $E_i \subset X_i$ so that $E_i \to E$ in the Gromov topology.  The subgroups $h_i(\text{Stab}_{G}(E))$ are abelian, and fix the flat $E_i$, for sufficiently large $i$. Thus $h_i(\text{Stab}_g(E))$ is contained in a unique maximal abelian subgroup $A_{E_i}$ of $\Gamma$.  If we fix a finite subset $W$ of $\text{Stab}_{G}(E)$ and $\epsilon > 0$, then for sufficiently large $i$, there is an automorphism $\sigma_i  \co A_{E_i} \to A_{E_i}$ so that
\begin{enumerate}
\item For every $w \in W$, and every $r_i \in E_i$,
\[	d_{X_i}(r_i, h_i(\sigma_i(w)).r_i) < \epsilon	.	\]
\item For any $k \in \text{Stab}_G(E)$ which acts trivially on $E$ we have $\sigma_i(h_i(k)) = h_i(k)$.
\end{enumerate}
The proof of the existence of such a $\sigma_i$ is the same as the proof of \cite[Proposition 5.4]{CWIF2}.  
Such a $\sigma_i$ induces a move of Type (B1) in a straightforward manner, since the adjacent edge groups to the vertex group $\text{Stab}_G(E)$ contain elements which act trivially on $E$, therefore we replace $h_i$ by the homomorphism which agrees with $h_i$ on all edge groups and on all vertex groups which are not $\text{Stab}_G(E)$, and replaces $h_i|_{\text{Stab}_G(E)}$ by $\sigma_i \circ h_i|_{\text{Stab}_G(E)}$.

We now construct shortening elements for all but finitely many of the intervals $[y_n, h_n(u).y_n]$ by following the proof of \cite[Theorem 5.1]{RipsSelaGAFA} (see \cite[Section 5]{CWIF2} for more details).

Finally, we are left with the case where $[y,u.y]$ is contained entirely in the discrete part of $T$.  We follow the proof of \cite[Theorem 6.1]{CWIF2}, which in turn followed \cite[Section 6]{RipsSelaGAFA}.  This argument naturally splits into a number of cases.

\medskip
{\bf Case 1}\qua $y$ is contained in the interior of an edge $e$.

{\bf Case 1a}\qua  $e$ is not contained entirely in a line $p_E \in \P$ and $\bar{e} \in T/G$ is a splitting edge.  

This case follows directly as in \cite[Section 6]{CWIF2}.

{\bf Case 1b}\qua $e$ is not completely contained in a line $p_E \in \P$ and $\bar{e}$ is not a splitting edge.

This case also follows directly as in \cite[Section 6]{CWIF2}.

{\bf Case 1c}\qua $e$ is contained in a line $p_E \in \P$ and $\bar{e}$ is a splitting edge.  

In this case, we have a graph of groups decomposition $G = H_1 \ast_{A_E} H_2$, where $A_E = \text{Stab}_G(E)$ (and $E$ is the flat in $\Cinf$ which projects to $p_E$).  The Dehn twist which is found in \cite[Section 6]{CWIF2} is naturally replaced by a bending move of type (B2).

{\bf Case 1d}\qua $e$ is contained in a line $p_E \in \P$ and $\bar{e}$ is not a splitting edge.

Once again the Dehn twist is replaced by a bending move of Type (B2).

\medskip
{\bf Case 2}\qua $y$ is a vertex of $T$.

Once again here there are four cases, depending on whether on edge adjacent to $y$ is or is not a splitting edge and is or is not contained in a line $p_E \in \P$.  In case an edge $e$ is not contained in a line $p_E \in \P$, we proceed exactly as in \cite{CWIF2}, following 
\cite{RipsSelaGAFA} directly.  In case $e \subset p_E$, we replace the shortening Dehn twists by bending moves of type (B2) as in Case 1 above.

Therefore, in any case, we can find moves which shorten all but finitely many of the $h_i$, as required.
\end{proof}

\section{Shortening quotients} \label{ShortenQSection}

In this section we recall the concept of {\em shortening quotients} from \cite[Section  5]{Sela1} and \cite[Definition 1.14]{SelaHyp}, and generalise this notion to the relatively hyperbolic setting.  The use of bending
moves means that the output of this section is slightly different to that 
in the hyperbolic case (though it will be good enough for our purposes).

Let $G$ be a finitely generated group, $\Gamma$ a torsion-free relatively hyperbolic group with abelian parabolics and $\{ h_n \co G \to \Gamma \}$ a stable sequence of homomorphisms, with associated $\Gamma$--limit group $L$.\footnote{{\em Stable} sequences of homomorphisms were defined in Definition \ref{Stable}.}  Suppose
that $L$ is nonabelian and freely indecomposable.  Let $\pi \co G \to 
L$ be the canonical quotient map.  The shortening procedure (given
below) constructs a sequence of homomorphisms $\{ \nu_i \co F 
\to \Gamma \}$, where $F$ is a finitely generated free group.  The 
sequence $\{ \nu_i \}$ has 
a subsequence converging to a $\Gamma$--limit group $Q$, which
 comes equipped with a canonical epimorphism
$\eta : L \to Q$.  See Lemma \ref{LtoQ} and Proposition 
\ref{ShortQProps} for more information.

We follow the construction from \cite[Section 3]{Sela1} and \cite[Section 3]{SelaHopf} (see \cite[Section  7]{CWIF2} for more details in the current context).  See also \cite{Serre} for more information on the Bass--Serre
theory, in particular \cite[I.5]{Serre} for the standard presentation
of the fundamental group of a graph of groups.

Given the situation described in the previous paragraph, we now describe the construction of $\{ \nu_i \}$, $Q$ and $\eta$.  Let 
$\Lambda_{L}$ be an abelian JSJ decomposition for $L$, with vertex groups $V^1, \ldots , V^m$ and edge groups $E^1, \ldots , E^s$.  Let $t^1, \ldots , t^b$ be the Bass--Serre generators for $L$ with respect to some (fixed) maximal subtree of $\Lambda_{L}$.
Note that $L$ is finitely generated, so each of the vertex groups $V_i$
is generated by a finite set together with the adjacent edge groups.

As in \cite[Theorem 3.2, comment at the bottom of page 45]{Sela1} we have not yet proved that the edge groups are finitely generated (though this will eventually turn out to be the case; see Corollary \ref{Abelian-fg} below). 

We can `approximate' the  finitely generated group $L$ by finitely presented groups $U_n$, each equipped with a graph of groups decomposition $\Lambda_n$ which is a `lift' of $\Lambda_{L}$.  We describe this approximation now, following \cite[Theorem 3.2]{Sela1} and \cite[Theorem 7.9]{CWIF2}.

Let $g_1, \ldots , g_k$ generate $G$ and let $\pi \co G \to L$ be the canonical quotient map associated to the convergent subsequence of $\{ h_n \co G \to \Gamma \}$.  Note that $L$ is clearly generated by $\{ \pi(g_1), \ldots , \pi(g_k) \}$. Now choose elements 
\[	v_1^1, \ldots , v_{l_1}^1, \ldots , v_1^m, \ldots  , v_{l_m}^m, t_1, \ldots , t_b \in G ,	\]
so that (i) $\pi(v_k^i) \in V^i$; (ii) $\pi(t_k) = t_k$; and (iii) for each generator $g_j$ of $G$ there is a word:
\[	g_j = w_j(v_1^1, \ldots , v_{l_m}^m, t_1, \ldots , t_b)	.	\]
Thus $G$ is generated by $\{ v_1^1, \ldots , v_{l_m}^m, t_1, \ldots , t_b \}$.  We may also assume, by adding finitely many more elements to our list if necessary, that each $V^i$ is generated by $\{ \pi(v_1^i), \ldots , \pi(v_{l_i}^i) \}$, together with the edge groups $E^j$ which are adjacent to $V^i$.  In case the edge groups adjacent to $V^i$ are all finitely generated, we assume that $V^i$ is in fact generated by $\{ \pi(v_1^i), \ldots , \pi(v_{l_i}^i) \}$.  For each $j$, let $e_1^j, e_2^j, \ldots \in G$ be a set of elements for which $\pi(e_p^j) \in E^j$ so that $E^j$ is generated by $\{ \pi(e_1^j), \pi(e_2^j), \ldots \}$.

We now define the groups $U_n$.  First, define the group $H_n$ to be the group with the generating set:
\[	\{ x_1^1, \ldots , x_{l_1}^1, \ldots , x_1^m, \ldots , x_{l_m}^m, y_1, \ldots , y_b, z_1^1, \ldots , z_1^s, z_n^1, \ldots , z_n^s \}	,	\]
together with the relations $[ z_{p_1}^j, z_{p_2}^j] = 1$ for $j = 1, \ldots , s$, and $1 \le p_1, p_2 \le n$.  Clearly there exists a natural epimorphism $\sigma_n \co H_n \to L$ defined by $\sigma_n(x_p^i) = \pi(v_p^i), \sigma_n(y_r) = t^r$ and $\sigma_n(z_d^j) = \pi(e_d^j)$.  We define $U_n$ to be the quotient of $H_n$ defined by adding as relations each word $w$ in the given generators of $H_n$ for which
\begin{enumerate}
\item $\sigma_n(w) = 1$,
\item the length of $w$ in the given generating set for $H_n$ is at most $n$ and
\item for some fixed index $i \in \{ 1 ,\ldots , m \}$, the word $w$ is a word in
\begin{enumerate}
\item the generators $x_1^i, \ldots, x_{l_i}^i$,
\item the elements $z_1^j, \ldots , z_n^j$ for any of the indices $j \in \{ 1 ,\ldots , s \}$ for which $E^j < V^i$ and
\item the words $y_rz_1^jy_r^{-1}, \ldots , y_rz_n^j y_r^{-1}$ for any pair of indices $(j,r)$ for which $t^rE^j (t^r)^{-1} < V^i$.
\end{enumerate}
\end{enumerate}

There exists a natural map $\kappa_n \co U_n \to U_{n+1}$.  Since $L$ is the quotient of $G$ by the stable kernel of $\{ h_i \co G \to \Gamma \}$, for each $n$ there exists $k_n > k_{n-1}$ so that the homomorphism $h_{k_n} \co G \to \Gamma$ induces a homomorphism $\lambda_n \co U_n \to \Gamma$ defined by $\lambda_n(x_p^i) = h_{k_n}(v_p^i)$, $\lambda_n(y_r) = h_{k_n}(t_r)$ and $\lambda_n(z_d^j) = h_{k_n}(e_d^j)$.  Also, the homomorphism $\sigma_n \co H_n \to L$ factors through $U_n$, and we denote the associated homomorphism from $U_n$ to $L$ by $\sigma_n$.

Since the second set of defining relations of $U_n$ consists of words whose letters are mapped by $\sigma_n$ into the same vertex group of $\Lambda_{L}$, each of the groups $U_n$ admits an abelian splitting $\Lambda_n$ which projects by $\sigma_n$ into the abelian decomposition $\Lambda_{L}$ of $L$.  That is to say each of the vertex groups $V_n^i$ in $\Lambda_n$ satisfies $\sigma_n(V_n^i) \le V^i$, each of the edge groups $E_n^j = \langle z_1^j, \ldots , z_n^j \rangle$ satisfies $\sigma_n(E_n^j) \le E^j$ and each of the Bass--Serre generators in $\Lambda_n$ satisfies $\sigma_n(y_r) = t^r$.

By the {\em surface} vertex groups of $\Lambda_n$ we mean the
vertex groups in $\Lambda_n$ which project to surface vertex groups
of $\Lambda_L$.  Since such groups are finitely presented, for all
but finitely many $n$ each surface vertex group in $\Lambda_n$ is
isomorphic to the corresponding surface vertex group of $\Lambda_L$,
via the natural projection from $U_n$ to $L$.
Denote by $\text{Mod}(U_n)$ the subgroup of $\text{Aut}(U_n)$ generated by
\begin{enumerate}
\item inner automorphisms,
\item Dehn twists in edge groups of $\Lambda_n$,
\item Dehn twists in cyclic groups which arise as the edge group
in a splitting of $U_n$ obtained by cutting the surface associated
to a surface vertex group of $\Lambda_n$ along a weakly essential
s.c.c.\footnote{Following \cite{RipsSelaJSJ}, a {\em weakly essential s.c.c} on a surface is a simple closed curve which is not null-homotopic, not boundary parallel and not the core of a M\"obius band.} and
\item generalised Dehn twists induced by the decomposition 
$\Lambda_n$.
\end{enumerate}

Let $W_n$ be the subgroup of $U_n$ generated by the $x_p^i$'s and the $y_r$'s.  Clearly the homomorphism $\kappa_n \co U_n \to U_{n+1}$ restricts to an epimorphism $\kappa_n \co W_n \to W_{n+1}$, and $\lambda_n$ restricts to a homomorphism from $W_n$ to $\Gamma$.

Thus we have the following commutative diagram, where $\iota_n \co W_n \to U_n$ is inclusion.
\[
\begin{CD}
W_1		@>\kappa_1>> 	W_2 @. \  \cdots 	@>\kappa_{n-2}>> 	W_{n-1} 	@>\kappa_{n-1}>> 	W_n \\
@V\iota_1VV		@V\iota_2VV	@.			@V\iota_{n-1}VV	@V\iota_nVV \\
U_1		@>\kappa_1>> 	U_2 @. \  \cdots 	@>\kappa_{n-2}>> 	U_{n-1} 	@>\kappa_{n-1}>> 	U_n \\
@V\lambda_1VV		@V\lambda_2VV	@.			@V\lambda_{n-1}VV	@V\lambda_nVV \\
\Gamma		@.		\Gamma	@. 	@. \Gamma		@. \Gamma
\end{CD}
\]
The $\Gamma$--limit group $L$ is the direct limit of the sequence 
$\{ (W_i, \kappa_i) \}$.  Denote the free group $W_1$ by $F$, and let 
$\phi \co F \to L$ be the canonical quotient map.  

Define an equivalence relation on $\text{Hom}(W_n,\Gamma)$ analogously to Definition \ref{ShortHomo}, using (i) automorphisms $\phi \in \text{Mod}(U_n)$; (ii) bending moves (where we restrict attention
to bending moves defined using the decomposition $\Lambda_n$ of $U_n$); and (iii) conjugations in $\Gamma$.  Also, define the norm$\| h \|$ of $h \in \text{Hom}(W_n,\Gamma)$, using the generating set $\{ x_1^1, \ldots , x_{l_m}^m,y_1,\ldots , y_r \}$.  For each $n \ge 1$, choose a homomorphism $\widehat{\lambda_n} \co W_n \to \Gamma$
so that (i) $\widehat{\lambda_n} \sim \lambda_n|_{W_n}$; and (ii) 
$\| \widehat{\lambda_n} \|$ is minimal amongst all homomorphisms equivalent to $\lambda_n|_{W_n}$.

Let $\xi_n$ be the given
epimorphism from $F$ to $W_n$ (obtained by identifying their generating sets in the obvious way), and let $\nu_n = \widehat{\lambda_n} \circ \xi_n \co F \to \Gamma$.  Passing to a subsequence of $\{ \nu_n \}$, we obtain an associated $\Gamma$--limit group $Q$, with
canonical quotient map $\nu_\infty \co F \to Q$.

\begin{lemma} \label{LtoQ}
The natural map between the generating set of $L$ and the generating set of $Q$ extends to an epimorphism $\eta : L \to Q$.  Furthermore,
$\eta$ is injective on each rigid vertex group of $\Lambda_L$.
\end{lemma}
\begin{proof}
The group $L$ is the direct limit of the sequence $\{ (W_i, \kappa_i) \}$, so the kernel of the map from $F$ to $L$ is $\SK(\xi_n)$. 

The $\Gamma$--limit group $Q$ is obtained from the sequence $\{ \nu_n = \widehat{\lambda_n} \circ \xi_n \}$.  It is clear that 
$\SK(\nu_n) \subseteq \SK(\xi_n)$.  This proves that there is an epimorphism $\eta \co L \to Q$ as required by the statement of the lemma.

Finally, we prove that the rigid vertex groups of $\Lambda_L$ are
mapped injectively into $Q$ by $\eta$.  First note that
rigid vertex groups are non-abelian, by definition.
Let $V^i$ be a rigid vertex
group of $L$, and let $x,y$ be distinct elements of $V^i$.  Suppose
that $u, v \in F$ are such that $\phi(u) = x$, $\phi(v) = y$.  By the properties of the above construction, there
exists $j_0 \ge 1$ so that for all $j \ge j_0$ the following hold: 

(i)\qua
$\xi_j(u), \xi_j(v) \in V_j^i$;\qua (ii)\qua $\xi_j(u) \neq \xi_j(v)$;\qua and 
\qua(iii)\qua
$\lambda_j(\xi_j(u)) \neq \lambda_j(\xi_j(v))$.

Now, the effect of all Dehn twists and generalised Dehn twists arising
from $\Lambda_j$, and also of the bending moves arising from 
$\Lambda_j$, is conjugation on $V_j^i$ (by the same element for 
all of $V_j^i$).  Hence $\widehat{\lambda_j}(\xi_j(u)) \neq
\widehat{\lambda_j}(\xi_j(v))$, which is to say that $\nu_j(u) \neq
\nu_j(v)$.  Since this is true for all $j \ge j_0$, we see that
$\nu_\infty(u) \neq \nu_\infty(v)$.  By construction, $\eta(x) =
\nu_\infty(u)$ and $\eta(y) = \nu_\infty(v)$.  Therefore,
$\eta$ is injective on $V^i$, as required.
\end{proof}

The group $Q$ is called the {\em shortening quotient} of $L$ associated
to $\{ h_n \co G \to \Gamma \}$.

Although we speak of {\em the} shortening quotient, $Q$ and $\eta$ depend on the choices of shortest homomorphism in the equivalence class of $\lambda_n$, and also on the convergent subsequence of 
$\{ \nu_n \}$ chosen. 

We now record some important properties of shortening quotients
for use later in this paper. 

\begin{proposition} \label{ShortQProps}
Suppose that $G$ is a finitely generated group, that $\Gamma$
is a torsion-free relatively hyperbolic group with abelian parabolics,
and that $\{ h_n \co G \to \Gamma \}$ is a sequence of homomorphisms
that converges into a faithful action of a non-abelian, freely
indecomposable strict $\Gamma$--limit group $L$ on an $\R$--tree.

Suppose that $\{ \nu_{n_t}  \co F \to \Gamma \}$ is a (convergent)
sequence of homomorphisms so that $\nu_{n_t} = 
\widehat{\lambda_{n_t}}
\circ \xi_{n_t}$, where $\widehat{\lambda_{n_t}} \co W_{n_t} \to 
\Gamma$ is a 
short homomorphism, as in the construction of shortening quotients
above.  Let $Q$ be the resulting shortening quotient, with canonical epimorphism $\eta \co L \to Q$ as in Lemma \ref{LtoQ} above.
\begin{enumerate}
\item  \label{ProperQ} If the sequence $\{ \nu_{n_t} \}$ is not contained 
in finitely many
conjugacy classes then $\eta$ is not injective, so $Q$ is a proper
quotient of $L$.
\item \label{Factor} If the sequence $\{ \nu_{n_t} \}$ is contained in finitely many
conjugacy classes, then there is some subsequence $\{ h_r \}$ of
$\{ h_{n_t} \}$ so that each $h_r$ factors through the canonical
epimorphism $\pi \co G \to L$.
\end{enumerate}
\end{proposition}
\begin{proof}
First we prove (\ref{ProperQ}). Suppose that $Q$ is not a proper
quotient of $L$, and that $\{ \nu_{n_t} \}$ is not contained in finitely
many conjugacy classes.  Then the limiting action of $Q$ on the
$\R$--tree $T$, induces a faithful action of $L$ on $T$.  This, in turn,
induces an abelian splitting $\Theta$ of $L$.  We use the splitting
$\Theta$ to shorten all but finitely many of the homomorphisms
$\widehat{\lambda_{n_t}}$, which contradicts the fact that they were
chosen to be short.  We argue as in Section \ref{ShortenSection}.

Consider a generator $u$ of $L$, the basepoint $y$ of $T$ and the segment $[y,u.y]$.  If $[y,u.y]$ intersects an IET component of $T$ in
more than a point, then
 \cite[Theorem
5.1]{RipsSelaGAFA} finds a
mapping class of the associated surface in the splitting of $L$ which
shortens the limiting action on $T$.  By \cite[Corollary 4.4]{CWIF2},
the length of $[y,u.y]$ in $T$ is the same as its the pre-image in 
$\mathcal{C}_\infty$, so we use the same automorphism to shorten
the limiting action on $\mathcal{C}_\infty$.
The surface group found above is conjugate into a surface
vertex group of the JSJ decomposition of $L$, and so the shortening
automorphism is a product of Dehn twists in s.c.c. on this surface.
Since a surface vertex is finitely presented, all but finitely many of the
$U_n$ have an isomorphic copy of this surface vertex, and we
can shorten the approximating homomorphisms.  Thus we
may suppose that $[y,u.y]$ does not intersect any IET components
in more than a point.

If $[y,u.y]$ intersects an axial subtree in a segment of positive length
then we apply the shortening moves in Section \ref{ShortenSection}
to shorten the limiting action on $\mathcal{C}_\infty$.  Once again,
these shortening moves can easily be lifted to all but finitely many
of the $U_n$ (by using the fact that the splitting $\Lambda_L$ of $L$
essentially encodes all abelian splittings of $L$).  Thus the 
appropriate shortening moves will shorten all but finitely many of the
approximating homomorphisms from $W_n$ to $\Gamma$.  It is worth
noting that even though we do not know that abelian subgroups of $L$
are finitely generated, the generalised Dehn twists that we need can be
chosen to be supported on a finitely generated subgroup, which is
finitely presented and will lift to all but finitely many of the $U_n$.

Finally, there is the case that $[y,u.y]$ is entirely contained in 
the discrete part of $T$.  In this case, the limiting action is not
shortened, but the appropriate Dehn twist, generalised Dehn
twist or bending move can be lifted to all but finitely many of
the $U_n$, and used to shorten the approximating homomorphisms,
just as in Section \ref{ShortenSection}.  The bending moves in any
of the cases in the proof are naturally interpreted using the 
homomorphisms from $W_n$ to $\Gamma$ (and the splitting
$\Lambda_n$ of $U_n$).  Thus (\ref{ProperQ}) is proved by
following the proof of Theorem \ref{Short}.

We now prove (\ref{Factor}).
Suppose that the elements of the sequence $\{ \nu_{n_t} \}$ belong to finitely
many conjugacy classes.  Since each of the $\widehat{\lambda}_{n_t}$
is short, the sequence $\| \nu_{n_t} \|$ is bounded, and by passing 
to a subsequence we may assume that $\{ \nu_{n_t} \}$ 
is constant, and equal to $\nu$, say.  For ease of notation, we refer 
to the constant subsequence of $\{ \nu_{n_t} \}$ by $\{ \nu_i \}$, and 
the associated subsequence of $\{ h_{n_t} \}$ by $\{ h_i \}$.

We have $Q \cong \nu(F) \le \Gamma$.  By Lemma \ref{LtoQ}, each of the rigid vertex groups of $\Lambda_L$ embed in $Q$, and hence in
$\Gamma$.  Let $V_j$ be a rigid vertex group in $\Lambda_L$, and 
denote by $V_j'$ the isomorphic image of $V_j$ in $\Gamma$ given by the above embedding.  By Lemma \ref{Abfg}, the edge groups in
$\Lambda_L$ adjacent to rigid vertex groups are finitely generated.
The edge groups adjacent to surface groups are finitely generated,
and $\Lambda_L$ does not have a pair of adjacent abelian vertex groups.  Therefore, all edge groups in $\Lambda_L$ are finitely generated, and therefore all vertex groups are also finitely generated.

Since all of the edge groups of $\Lambda_L$ are finitely generated, we assumed during the construction of the $U_i$ and $W_i$ that the vertex group $V_j$ of $\Lambda_L$ was generated by
$\{ \pi(v_1^j), \ldots , \pi(v_{l_j}^j) \}$, for each $1 \le j \le m$. 
Now, for each $1 \le j \le m$ and each $1 \le k \le l_j$, we have
$h_i(v_k^j) = \lambda_i(x_k^j)$.  Therefore, if $Y_j = \langle v_1^j, \ldots , v_{l_j}^j \rangle \le G$ and $Z_j^i = \langle x_1^j, \ldots , x_{l_j}^j \rangle \le W_i$ then
\[	h_i(Y_j) = (\lambda_i \circ \iota_i) (Z_j^i)	.	\]
By construction, $(\lambda_i \circ \iota_i) \sim \widehat{\lambda}_i$.

The fact that $Q \cong \nu(F) \le \Gamma$ implies that 
$\widehat{\lambda}_i(Z_j^i)$ is isomorphic to $V_j'$.  Now,
$(\lambda_i \circ \iota_i)$ and $\widehat{\lambda}_i$ are related by 
a sequence of
shortening moves, each of which acts by conjugation on $Z_j^i$.
Thus $(\lambda_i \circ \iota_i) (Z_j^i)$ is conjugate to $V_j'$ in 
$\Gamma$. The same is therefore true of $h_i(Y_j)$.

Since the edge groups of $\Lambda$ are finitely generated (and abelian), and since the non-rigid vertex groups are finitely presented, a presentation for $L$
can be found which consists of presentations for the rigid vertex 
groups, and finitely many other generators and relations.  The 
sequence $\{ h_n \co G \to \Gamma \}$ converges to $\pi \co G 
\to L$, so each of these finitely many extra relations will be sent 
to $1$ by $h_l$ for sufficiently large $l$.  However,
we have already seen that the relations corresponding to the
rigid vertex groups are mapped to $1$ by all the $h_{i}$.
This shows that all but finitely many of the $h_i$ factor 
through $\pi \co G \to L$, as required.  This proves
(\ref{Factor}).
\end{proof}

\section{$\Gamma$--limit groups} \label{MainSection}

In this section we follow \cite[Section  1]{SelaHyp} in order to understand $\Gamma$--limit groups, and to understand $\text{Hom}(G,\Gamma)$, where $G$ is an arbitrary finitely generated group.  The main technical results of this section are Theorem \ref{SeriesTerminates} and Theorem \ref{GeneralFiniteRes}.

These technical results are then applied to yield various applications:  Theorem \ref{FullyResiduallyGamma}, Proposition \ref{CountablyMany}, Theorem \ref{FiniteSubsystem}, and in the next section the construction of Makanin--Razborov diagrams.

\begin{remark}
A previous version of this paper claimed that any strict $\Gamma$--limit group which is freely indecomposable and nonabelian admits a principal cyclic splitting.  However, this is false.

To see this, consider the following example.  Let $\Gamma$ be the fundamental group of a finite volume hyperbolic $n$--manifold $M$ with torus cusps, for $n \ge 3$, so that $\Gamma$ admits no splitting over $\mathbb Z$.  Let $M'$ be the double of $M$ along one of its cusps, and $\Gamma' = \pi_1(M')$.  By repeatedly Dehn twisting along a simple closed curve in the attaching torus, it is not difficult to see that $\Gamma'$ is a strict $\Gamma$--limit group.  However, the JSJ decomposition of $\Gamma'$ is given by two copies of $\Gamma$ attached by a (noncyclic) abelian edge group.  It is clear that $\Gamma'$ admits no principal cyclic splitting.

It is possible to adapt the proof of \cite[Theorem 3.2]{Sela1} to this context to prove that the above example is one of the few ways in which a strict $\Gamma$--limit group can fail to have a principal cyclic splitting.  However, rather than make this long diversion, we prefer to work entirely with the abelian JSJ decomposition of a strict $\Gamma$--limit group $L$, which we already know to be nontrivial, so long as $L$ is freely
indecomposable, non-abelian and not a surface group, by Proposition
\ref{JSJnontriv}.
\end{remark}

Let $G$ be a fixed finitely generated group.  Define an order on the set of $\Gamma$--limit groups that are quotients of $G$ as follows:  suppose $R_1$ and $R_2$ are both $\Gamma$--limit groups that are quotients of $G$, and that $\eta_i \co G \to R_i$ are the (fixed) canonical quotient maps.  We say $R_1 > R_2$ if there exists an epimorphism with non-trivial kernel $\tau \co R_1 \to R_2$ so that $\eta_2 = \tau \circ \eta_1$.  We say that $R_1$ and $R_2$ are {\em equivalent} if there is an isomorphism $\tau \co R_1 \to R_2$ so that $\eta_2 = \tau \circ \eta_1$.

The following is one of the main technical results of this paper.

\begin{theorem}[Compare {\cite[Theorem 1.12]{SelaHyp}}]  
\label{SeriesTerminates}
Let $\Gamma$ be a torsion-free group which is hyperbolic relative to free abelian subgroups.  Every decreasing sequence of $\Gamma$--limit groups:
\[	R_1 > R_2 > R_3 > \ldots	,	\]
which are all quotients of a finitely generated group $G$,
terminates after finitely many steps.
\end{theorem}

For limit groups, the analogous result has a short proof using the fact that free groups are linear (see \cite[Corollary 1.9, page 4]{BFSela}).  However (as observed by M. Kapovich; see, for example, \cite[Subsection 1.4]{BestvinaQuestions}), not all hyperbolic groups are linear and therefore not all relatively hyperbolic groups are linear. 

It is worth noting that we do not insist that the $\Gamma$--limit groups $R_n$ in the statement of Theorem \ref{SeriesTerminates} be strict $\Gamma$--limit groups.  This will be important later, because to study a single homomorphism, we consider a constant sequence, which leads to a $\Gamma$--limit group which need not be strict.

Before we prove Theorem \ref{SeriesTerminates}, we prove the following lemma (implicit in \cite[page 7]{SelaHyp}):

\begin{lemma} \label{d-generated}
Let $\Xi$ be a finitely generated group, let $L$ be a $\Xi$--limit group and suppose that $L$ is $d$--generated.  Then $L$ can be obtained as a limit of homomorphisms $\{ f_n \co F_d \to \Xi \}$, where $F_d$ is the free group of rank $d$.
\end{lemma}
\begin{proof}
The group $L$ is obtained as the limit of a (convergent)
sequence of homomorphisms
 $\{ h_n \co G \to \Xi \}$, where $G$ is a finitely generated group.  Let
$\eta \co G \to L$ be the canonical quotient map. Let
$\{ r_1, \ldots , r_d \}$ be a generating set for $R$ and let
$\{ y_1 , \ldots , y_d \} \subset G$ be such that $\eta(y_i) = r_i$
for $1 \le i \le d$.  Let $G_d \le G$ be the subgroup generated by
$\{ y_1 , \ldots , y_d \}$, and let $\pi \co F_d \to G_d$ be the
natural quotient map.  Define $h_n = f_n \circ \pi$.  It is 
easy to see that the sequence $\{ h_n \co F_d \to \Xi \}$ converges
to $L$.
\end{proof}

\begin{proof}[Proof of Theorem \ref{SeriesTerminates}]
We follow the outline of the proof of \cite[Theorem 1.12]{SelaHyp}.  In order to obtain a contradiction, we suppose that there exists a finitely generated group $G$ and an infinite descending sequence of $\Gamma$--limit groups:
\[	R_1 > R_2 > R_2 > \ldots	,	\]
all quotients of $G$.

Without loss of generality, we may assume $G = F_d$, the free group of rank $d$.  Let $\{ f_1, \ldots, f_d \}$ be a basis for $F_d$, and let $C$ be the Cayley graph of $F_d$ with respect to this generating set.  We construct a particular decreasing sequence of $\Gamma$--limit groups as follows.  Let $T_1$ be a $\Gamma$--limit group with the following properties:
\begin{enumerate}
\item $T_1$ is a proper quotient of $F_d$.
\item $T_1$ can be extended to an infinite decreasing sequence of $\Gamma$--limit groups: $T_1 > L_2 > L_3 > \ldots$.
\item The map $\eta_1 \co F_d \to T_1$ maps to the identity the maximal number of elements in the ball of radius $1$ about the identity in $C$ among all possible epimorphisms from $F_d$ to a $\Gamma$--limit group $L$ that satisfies the first two conditions.
\end{enumerate}
Continue to define the sequence inductively.  Suppose that the finite sequence $T_1 > T_2 > \ldots > T_{n-1}$ has been constructed, and choose $T_n$ to satisfy:
\begin{enumerate}
\item $T_n$ is a proper quotient of $T_{n-1}$.
\item The finite decreasing sequence of $\Gamma$--limit groups $T_1 > T_2 > \ldots > T_n$ can be extended to an infinite decreasing sequence.
\item The map $\eta_n \co F_d \to T_n$ maps to the identity the maximal number of elements in the ball of radius $n$ about the identity in $C$ among all possible maps from $F_d$ to a $\Gamma$--limit group $L_n$ satisfying the first two conditions.
\end{enumerate}

Since each of the $\Gamma$--limit groups $T_n$ is a quotient of $F_d$, each $T_n$ is $d$--generated. Let $\{ r_{1,n}, \ldots , r_{d,n} \}$ be a generating set for $T_n$.  By Lemma \ref{d-generated}, $T_n$ can be obtained as a limit of a sequence of homomorphisms $\{ v_i^n \co F_d \to \Gamma \}$, with the quotient map $\eta_n \co F_d \to T_n$ sending $f_i$ to $r_{i,n}$.

For each $n$, choose a homomorphism $v_{i_n}^n \co F_d \to \Gamma$ for which:
\begin{enumerate}
\item Every element in the ball of radius $n$ about the identity in $C$ that is mapped to the identity by $\eta_n \co F_d \to T_n$ is mapped to the identity by $v_{i_n}^n$.  Every such element that is mapped to a nontrivial element by $\eta_n$ is mapped to a nontrivial element by $v_{i_n}^n$.
\item There exists an element $f \in F_d$ that is mapped to the identity by $\eta_{n+1} \co F_d \to T_{n+1}$ for which $v_{i_n}^n(f) \neq 1$.
\end{enumerate}
Denote the homomorphism $v_{i_n}^n$ by $h_n$. By construction, the set of homomorphisms $\{ h_n \co F_d \to \Gamma \}$ does not belong to a finite set of conjugacy classes.  Therefore, from the sequence $\{ h_n \}$ we can extract a subsequence that converges to a (strict) $\Gamma$--limit group, denoted $T_\infty$.  By construction, the $\Gamma$--limit group $T_\infty$ is the direct limit of the sequence of (proper) epimorphisms:
\[	F_d \to T_1 \to T_2 \to \cdots\]
Let $\eta_\infty \co F_d \to T_\infty$ be the canonical quotient map.

The $\Gamma$--limit group $T_\infty$ is a proper quotient of each of the $\Gamma$--limit groups $T_n$.  For each index $n$, the group $T_n$ was chosen to maximise the number of elements in the ball of radius $n$ about the identity in $C$ that are mapped to the identity in $\Gamma$ among all $\Gamma$--limit groups that are proper quotients of $T_{n-1}$ and that admit an infinite descending chain of $\Gamma$--limit groups.  Therefore, it is not difficult to see that $T_\infty$ does not admit an infinite descending chain of $\Gamma$--limit groups.

The proof of the following proposition is similar to that of \cite[Proposition 1.16]{SelaHyp}.

\begin{proposition} \label{FiniteRes}
There is a subsequence $\{ h_{n_t} \}$ of $\{ h_n \co F_d \to \Gamma \}$ so that each $h_{n_t}$ factors through $\eta_\infty \co F_d \to T_\infty$.
\end{proposition}
\begin{proof}
Let $T_\infty = H_1 \ast \cdots \ast H_p \ast F$ be the Grushko
decomposition of $T_\infty$.  Consider $H_1$, and a 
finitely generated subgroup $F^1$ of $F_d$ so that $\eta(F^1) = H_1$.
Use the sequence $\{ h_n|_{F^1} \co F^1 \to \Gamma \}$ to construct
a shortening quotient of $H_1$, proper if possible.  
Repeat for each of the $H_i$ in 
turn, starting with the subsequence of $\{ h_n \}$ used to
construct the shortening quotient of $H_{i-1}$ as a starting point
for constructing a shortening quotient of $H_i$.
In this way, we eventually get a quotient $L_1$ of $T_\infty$ (map $F$ 
to itself in this quotient).

For each of the $H_i$ which admitted a proper
shortening quotient, take the Grushko decomposition of the 
corresponding free factor in $L_1$, and
for each of the non-free factors $H_j'$ in this decomposition, 
attempt to find a proper
shortening quotient of $H_j'$, passing to finer and finer subsequences
each time.

Continuing in this manner, we obtain a sequence of quotients:
\[T_\infty \to L_1 \to \cdots \to L_s\]
which are all proper quotients, except possibly in the case $s=1$ and
$T_\infty \to L_1$ is an isomorphism.
This process will terminate after finitely many steps, since $T_\infty$
does not admit an infinite descending chain of proper $\Gamma$--limit quotients. 

We now have a subsequence $\{ f_i \}$ of $\{ h_n \}$, the final
subsequence used in constructing a shortening quotient ($\{ f_i \}$ is
contained in all of the other subsequences used in this process). 
By construction, $L_s$ is a free product of finitely generated free
groups, and $\Gamma$--limit groups which do not admit a proper
shortening quotient starting from the sequence $\{ f_i \}$.

Let $L_s = K_1 \ast \ldots \ast K_q \ast F''$ be the given free
product decomposition of $L_s$, where the $K_i$ are not free,
and let $\sigma_i \co F_{k_i} \to K_i$ be the canonical quotient
map.
By Proposition \ref{ShortQProps}, there is a subsequence $\{f_i' \}$
of $\{ f_i \}$ so that all $f_i'|_{F_{k_1}} \co F_{k_1} \to
\Gamma$ factor through $\sigma_1$.  Similarly, there is a
further subsequence of $\{ f_i' \}$ which, when suitably restricted,
factor through $\sigma_2$.  By passing to further subsequences,
we can see that there is a subsequence of the $f_i$ which all
factor through $\eta_s \co F_d \to L_s$.

We have already noted that a
subsequence of the homomorphisms $\{ h_n \}$ can be assumed 
to factor through $L_s$. Thus, by induction on $j$, 
suppose that there is a further
subsequence of  $\{ h_n \}$ which factors through the canonical
quotient map $F_d \to L_{j+1}$.  The rigid vertex groups of the non-free factors of the Grushko decomposition of $L_j$ inject into $L_{j+1}$, by the shortening argument.  Since abelian subgroups of $\Gamma$ are finitely generated (Lemma \ref{Abfg}), we may assume by induction that the abelian subgroups of $L_{j+1}$ are finitely generated.  Thus the abelian vertex groups of the JSJ decomposition of a Grushko factor of $L_j$, as well as the edge groups, are finitely generated.  This implies that the abelian subgroups of $L_j$ are all finitely generated. The
canonical quotient map $L_j \to L_{j+1}$ embeds the rigid
vertex groups in the JSJ decomposition of $L_j$ into $L_{j+1}$.  Thus
we may assume that a subsequence of $\{ h_n \}$ factors through the
canonical quotient map from a finitely generated subgroup of
$F_d$ to the rigid vertex groups of $L_j$.
The other vertex groups in the JSJ decomposition of $L_j$ are finitely
presented, so there is some further subsequence of $\{ h_n \}$ which
factors through the canonical quotient to these vertex groups (the proof
of these claims are almost identical to the proof of Proposition 
\ref{ShortQProps}.(2)) .  It now
follows that some subsequence of $\{ h_n \}$ factors through the 
canonical quotient $F_d \to L_j$.
Thus, by induction, a subsequence of $\{ h_n \}$ factors through 
$\eta_\infty \co F_d \to T_\infty$, as required.
\end{proof}

The homomorphisms $\{ h_n \co F_d \to \Gamma \}$ were chosen so that for every index $n$ there exists some element $f \in F_d$ for which $\eta_{n+1}(f) = 1$ and $h_n(f) \neq 1$.  Now, $\eta_\infty \co F_d
\to T_\infty$ is the 
direct limit of the sequence 
\[F_d \to T_1 \to T_2 \to \cdots\]
where the induced map from $F_d$ to $T_i$ is $\eta_i$.  Thus, for
all $n \ge 1$ and all $f \in F_d$, if $\eta_{n+1}(f) = 1$ then
$\eta_\infty(f) = 1$.  By Proposition \ref{FiniteRes} it is possible to extract a subsequence $\{ h_{n_t} \co F_d \to \Gamma \}$ that factors through the map $\eta_\infty \co F_d \to T_\infty$, which is to say that there is a homomorphism $\pi_t \co T_\infty \to \Gamma$ so that $h_{n_t} = \pi_t \circ \eta_\infty$.  Hence, for every index $t$, and every element $f \in F_d$, if $\eta_{n_t+1}(f) = 1$ then $\eta_{\infty}(f) = 1$, which implies that $h_{n_t}(f) = 1$, in contradiction to the way that the homomorphisms $h_n$ were chosen.  This finally ends the proof of Theorem \ref{SeriesTerminates}.
\end{proof}

\begin{corollary} \label{LimitHopf}
Let $\Gamma$ be a torsion-free relatively hyperbolic group with abelian parabolics, and let $L$ be a $\Gamma$--limit group.  Then $L$ is Hopfian.
\end{corollary}
Corollary \ref{LimitHopf} generalises, and gives an independent
(though similar) proof of, \cite[Theorem A]{RelHyp}, and implies that the relation defined on $\Gamma$--limit groups which are quotients of a fixed group $G$ is a partial order. The proof of the following theorem is similar to that of \cite[Theorem 1.17]{SelaHyp}, and closely follows that of Proposition \ref{FiniteRes} above.

\begin{theorem} \label{GeneralFiniteRes}
Suppose that $G$ is a finitely generated group, and $\Gamma$ a torsion-free relatively hyperbolic group with abelian parabolics, and let
$\{ h_n \co G \to \Gamma \}$ be a sequence of homomorphisms converging to a $\Gamma$--limit group $L$, with canonical quotient $\pi \co G \to L$.  There 
exists a subsequence $\{ h_{n_t} \}$ of $\{ h_n \}$ so that each 
$h_{n_t} \co G \to \Gamma$ factors through $\pi \co G \to L$.
\end{theorem}
\begin{proof}
The sequence $\{ h_n \}$ converges to a faithful action of $L$ on an 
$\R$--tree $T$.  As in the proof of Theorem \ref{FiniteRes}, take the
Grushko decomposition of $L$, and attempt to form proper shortening
quotients of the non-free factors.  For each such proper quotient,
take the Grushko decomposition and attempt to form further
proper shortening quotients from the non-free factors, starting with
a finer subsequence of $\{ h_n \}$ at each stage.  We eventually 
obtain a sequence (finite by Theorem \ref{SeriesTerminates})
\[	G \to L \to L_2 \to \cdots \to L_s	.	\]
The rest of the proof is the same as that of Proposition \ref{FiniteRes}.
\end{proof}

Theorem \ref{GeneralFiniteRes} is the tool which we will use in the study
of $\Gamma$--limit groups as a replacement for the finite
presentability of limit groups.

\begin{corollary} \label{StartwithL}
Let $\Gamma$ be a torsion-free group which is hyperbolic relative to
abelian subgroups, and suppose that $L$ is a $\Gamma$--limit group.
Then there is a sequence of homomorphisms $\{ \rho_n \co L \to 
\Gamma \}$ which converge to $L$.
\end{corollary}
\begin{proof}
In case $L$ is not a strict $\Gamma$--limit group, it is isomorphic to
a finitely generated subgroup of $\Gamma$, so we can take a constant
sequence of embeddings.
In case $L$ is a strict $\Gamma$--limit group, the result follows from
Theorem \ref{GeneralFiniteRes}.
\end{proof}

Corollary \ref{StartwithL} simplifies the definition of shortening 
quotients considerably.  In particular, the groups $U_n$ and $W_n$
are no longer required.  We briefly describe the simplified construction
of shortening quotients.

Start with a non-abelian and freely
indecomposable strict $\Gamma$--limit group $L$, and a sequence
of homomorphisms $\{ h_n \co L \to \Gamma \}$ which converges to 
$L$.  For each $n$,
let $\widehat{h_n}$ be a short homomorphism which is equivalent
to $h_n$ (where now `short' is as defined in Definition \ref{ShortHomo}).
Passing to a convergent subsequence of $\{ \widehat{h_n} 
\co L \to \Gamma \}$, we find a $\Gamma$--limit group $Q$ which is
a quotient of $L$.

In this context, Proposition \ref{ShortQProps}.(1) follows immediately
from Theorem \ref{Short}.  The homomorphisms $h_n$ all have
$L$ as the domain, so Proposition \ref{ShortQProps}.(2) is
tautological in this context.

We record the simpler version of Proposition \ref{ShortQProps} here:

\begin{proposition} \label{ImprovedShortQ}
Suppose that $\Gamma$ is a torsion-free relatively hyperbolic group
with abelian parabolics, that $L$ is a non-abelian and freely-indecomposable strict $\Gamma$--limit group, and that $\{ h_n \co
L \to \Gamma \}$ is a sequence of homomorphisms converging to $L$.
Suppose further that $\widehat{h_{n_t}} \co L \to \Gamma$ is a convergent sequence so that $\widehat{h_{n_t}}$ is short and equivalent
to $h_{n_t}$, for each $t$.  Let $Q$ be the resulting shortening quotient,
with canonical epimorphism $\eta \co L \to Q$.

If the sequence $\widehat{h_{n_t}}$ is not contained in finitely many conjugacy classes then $Q$ is a proper quotient of $L$.
\end{proposition}

\begin{definition}
Let $\Gamma$ and $G$ be finitely generated groups.  We say that $G$ is {\em fully residually $\Gamma$} if for every finite $\mathcal F \subset G$ there is a homomorphism $h \co G \to \Gamma$ which is injective on $\mathcal F$.
\end{definition}

\begin{theorem} \label{FullyResiduallyGamma}
Let $\Gamma$ be a torsion-free relatively hyperbolic group with abelian parabolics.  A finitely generated group $G$ is a $\Gamma$--limit group if and only if it is fully residually $\Gamma$.
\end{theorem}
\begin{proof}
If $G$ is a fully residually $\Gamma$ then it is certainly a $\Gamma$--limit group. 
Conversely, suppose that $G$ is a $\Gamma$--limit group.  By Corollary 
\ref{StartwithL}, there is a sequence of homomorphisms
$\{ h_n \co G \to \Gamma \}$ which converges to $G$.  Using this 
sequence, it is clear that $G$ is fully residually $\Gamma$.
\end{proof}

\begin{proposition} \label{CountablyMany}
Let $\Gamma$ be a torsion-free relatively hyperbolic group with abelian parabolics.  Then there are only countably many $\Gamma$--limit groups.
\end{proposition}
\begin{proof}
If $\Gamma$ is abelian, then all $\Gamma$--limit groups are finitely generated abelian groups, of which there are only countably many.

Suppose then that $\Gamma$ is nonabelian.  Let $L$ be a $\Gamma$--limit group, and consider a sequence
of shortening quotients of $L$:
\[ L \to L_2 \to \ldots \to L_s	\]
constructed as in the proof of Theorem \ref{GeneralFiniteRes}, starting
with a sequence $\{ h_n \co L \to \Gamma \}$ which converges to $L$.
Note that none of the non-free factors of $L_s$ admits a proper 
shortening quotient, and if 
$s > 1$ then each map in the sequence is a proper quotient.  We will
prove by induction that there are only countably many 
$\Gamma$--limit groups.

For the base case, note that $L_s$ is a free product of a finitely
generated free group and a finite collection of finitely generated
subgroups of $\Gamma$.  There are only countably many such groups.

Assume by induction there are only countably
many groups which can be $L_{j+1}$ in the sequence.  Now, $L_j$ is a
free product of a finitely generated free group and freely indecomposable $\Gamma$--limit groups.  Some of the free factors
of $L_j$ embed into $L_{j+1}$ (if they do not admit a proper shortening
quotient).  If $H_k$ is one of the other non-free Grushko factors
of $L_j$, then $H_k$ admits an abelian JSJ decomposition 
$\Lambda_{H_k}$.
By Lemma \ref{LtoQ}, the rigid vertex groups of $\Lambda_{H_k}$
embed into $L_{j+1}$.  By induction and Lemma \ref{Abfg} we
may assume that abelian subgroups of $L_{j+1}$ are finitely generated.
By the argument in the proof of Proposition \ref{ShortQProps}, we see
that the edge groups and vertex groups of $\Lambda_{H_k}$ are finitely
generated.  Therefore each free factor of $L_j$ is either (i) a finitely
generated free group; (ii) a finitely generated subgroup of $L_{j+1}$; or
(iii) can be formed by
taking finitely many HNN extensions and amalgamated free products
of finitely generated subgroups of $L_{j+1}$ and finitely generated
free groups over finitely generated abelian
subgroups.  There are only countably many such constructions, and so
there are only countably many choices for $L_j$.
This also implies that abelian subgroups of $L_j$ are finitely generated
(which we assumed by induction for $L_{j+1}$).

Therefore, by induction, there are only countably many choices for 
$L$.
This completes the proof of the proposition.
\end{proof}

The following result was proved in the course of proving Proposition 
\ref{CountablyMany}.

\begin{corollary} \label{Abelian-fg}
Any abelian subgroup of a $\Gamma$--limit group is a finitely generated free abelian group.
\end{corollary}

\begin{proposition} [Compare Proposition 1.20, \cite{SelaHyp}] \label{MaxShorten}
Let $G$ be a finitely generated group and $\Gamma$ a torsion-free relatively hyperbolic group with abelian parabolics.  Let $R_1, R_2 \ldots $ be a sequence of $\Gamma$--limit groups that are all quotients of $G$ so that
\[	R_1 < R_2 < \cdots	.	\]
Then there exists a $\Gamma$--limit group $R$, a quotient of $G$, so that $R > R_m$ for all $m$.
\end{proposition}
\begin{proof}
For each $m$, choose a homomorphism $h_m \co G \to \Gamma$ that factors through the quotient map $\eta_m \co G \to R_m$ as $h_m = h'_m \circ \eta_m$ and so that $h_m'$ is injective on the ball of radius $m$ in $R_m$.  This is possible by Theorem \ref{GeneralFiniteRes}.

A subsequence of $\{ h_m \}$ converges to a $\Gamma$--limit group $R$, which is a quotient of $G$.  By Theorem \ref{GeneralFiniteRes}, we may assume that each element of this subsequence factors through the canonical quotient map $\eta \co G \to R$.

We prove that $R > R_m$ for each $m$.  We have quotient maps $\eta \co G \to R$, and $\eta_i \co G \to R_i$.  Since $R_i < R_{i+1}$ there exists $\tau_i \co R_{i+1} \to R_i$ so that $\eta_i = \tau_i \circ \eta_{i+1}$.  In particular, $\text{ker}(\eta_{i+1}) \subseteq \text{ker}(\eta_i)$.

Let $\A$ be the fixed finite generating set for $G$.  We attempt to define a homomorphism $\kappa_i \co R \to R_i$ as follows:  for $a \in \A$, define $\kappa_i(\eta(a)) = \eta_i(a)$.  This is well-defined if and only if $\text{ker}(\eta) \subseteq \text{ker}(\eta_i)$.  Therefore, suppose that $g \in \text{ker}(\eta)$.  Since each $h_j$ factors through $\eta$, we have $h_j(g) = 1$ for all $j$.  Suppose that $g$ lies in the ball of radius $n$ about the identity in the Cayley graph of $G$.  Then for all $j$, the element $\eta_j(g)$ lies in the ball of radius $n$ about the identity in $R_j$.  Since $h_j(g) = 1$, for all $j$, and by the defining property of the $h_j$, if $k \ge n$ then $\eta_k(g) = 1$.  Thus since for all $j$ we have $\text{ker}(\eta_{j+1}) \subseteq \text{ker}(\eta_j)$ we have $\eta_i(g) = 1$, as required.  We have constructed a homomorphism $\kappa_i \co R \to R_i$ so that $\eta_i = \kappa_i \circ \eta$, which is to say that $R > R_i$.  This finishes the proof.
\end{proof}

Propositions \ref{CountablyMany} and \ref{MaxShorten} imply that there are maximal elements for the set of $\Gamma$--limit groups which are quotients of a fixed finitely generated group $G$, under the order described before Theorem \ref{SeriesTerminates}.

Recall that we say that two $\Gamma$--limit groups which are quotients of $G$, $\eta_1 \co G \to R_1$ and $\eta_2 \co G \to R_2$ are {\em equivalent} if there is an isomorphism $\tau \co R_1 \to R_2$ so that $\eta_2 = \eta_1 \circ \tau$.

\begin{proposition} [Compare {\cite[Proposition 1.21]{SelaHyp}}] \label{FiniteMax}
Let $G$ be a finitely generated group and $\Gamma$ a torsion-free group hyperbolic relative to free abelian subgroups.    Then there are only finitely many equivalence classes of maximal elements in the set of $\Gamma$--limit groups that are quotients of $G$.
\end{proposition}
\begin{proof}
The following proof was explained to me by Zlil Sela in the context of torsion-free hyperbolic groups.  The same proof works in the current context.

Suppose on the contrary that there are infinitely many non-equivalent maximal $\Gamma$--limit groups $R_1, R_2, \ldots$, each a quotient of $G$.  Let $\eta_i \co G \to R_i$ be the canonical quotient map.  Fixing a finite generating set $\A$ for $G$, we fix a finite generating set for each of the $R_i$, and hence obtain maps $\nu_i \co F_d \to R_i$, where $d = | \A |$.  There is a fixed quotient map $\pi \co F_d \to G$ so that for each $i$ we have $\nu_i = \eta_i \circ \pi$.

For each $i$, consider the set of words of length $1$ in $F_d$ that are mapped to the identity by $\nu_i$.  This set is finite for each $i$, and there is a bound on its size, so there is a subsequence of the $R_i$ so that this set is the same for all $i$.  Starting with this subsequence, consider those words of length $2$ in $F_d$ which are mapped to the identity by $\nu_i$, and again there is a subsequence for which this (bounded) collection is the same for all $i$.  Continue with this process for all lengths of words in $F_d$, passing to finer and finer subsequences, and consider the diagonal subsequence.  We continue to denote this subsequence by $R_1, R_2, \ldots$.

Now, for each $i$, choose a homomorphism $h_i \co F_d \to \Gamma$ so that for words $w$ of length at most $i$ in $F_d$, we have $h_i(w) = 1$ if and only if $\nu_i(w) = 1$, and so that $h_i$ factors through the quotient map $\pi \co F_d \to G$.  This is possible because each $R_i$ is a $\Gamma$--limit group which is a quotient of $G$.

A subsequence of $\{ h_i \co F_d \to \Gamma \}$ converges into a $\Gamma$--limit group $M$, which is a quotient of $G$ since all $h_i$ factor through $\pi$.  Let $\psi \co F_d \to M$ be the canonical quotient, and $\phi \co G \to M$ the quotient for which $\psi = \phi \circ \pi$.  Note that a word $w$ of length at most $i$ in $F_d$ maps to the identity under $\psi$ if and only if $\nu_i(w) = 1$.

Now, $R_1, R_2, \ldots $ are non-equivalent maximal $\Gamma$--limit quotients, so (possibly discarding one $R_i$ which is equivalent to $M$) are all non-equivalent to $M$.  Therefore, for each $i$ there does not exist a homomorphism $\mu \co M \to R_i$ so that $\nu_i = \mu \circ \psi$.  That is to say that for each $i$ there exists $u_i \in F_d$ so that $\psi(u_i) = 1$ but $\nu_i(u_i) \neq 1$.

Let $\{ \tau_i \co F_d \to \Gamma \}$ be a sequence of homomorphisms that all factor through $\pi \co F_d \to G$ so that
\begin{itemize}
\item a word $w \in F_d$ of length at most $i$ satisfies $\tau_i(w) = 1$ if and only if $\nu_i(w) = 1$ and
\item $\tau_i(u_i) \neq 1$.
\end{itemize}
By Theorem \ref{FiniteRes} there is a subsequence $\{ \tau_{n_i} \}$ of $\{ \tau_i \}$ which converges into a $\Gamma$--limit group (which must be $M$) so that each $\tau_{n_i}$ factors through $\psi \co F_d \to M$.  Therefore, there is $r_i \co M \to \Gamma$ so that $\tau_{n_i} = r_i \circ \psi$.  However, we have that $\psi(u_{n_i}) = 1$, but $1 \neq \tau_{n_i}(u_{n_i}) = r_i(\psi(u_{n_i})) = r_i(1) = 1$, a contradiction.  This contradicts the existence of $R_1, R_2, \ldots$, and finishes the proof.
\end{proof}

One of the motivations for the work done in this paper and its
predecessors \cite{CWIF, CWIF2, RelHyp} were the following
questions asked by Sela \cite[I.8]{SelaProblems}:
Let $G$ be a CAT$(0)$ group with isolated flats (see \cite{Hruska}).

I.8(i)\qua Is $G$ Hopf?

I.8(ii)\qua Is it true that every system of equations over $G$ is equivalent to a finite system?

I.8(iii)\qua Is it possible to associate a Makanin--Razborov diagram to a system of equations over $G$?

I.8(iv)\qua Is it possible to decide if a system of equations over $G$ has a solution?

We believe that relatively hyperbolic groups with abelian parabolics
form a more natural context for this question than CAT$(0)$ groups
with isolated flats.  Note that `virtually abelian' rather than `abelian'
would be a natural class (containing CAT$(0)$ with isolated flats), 
but our methods only work for abelian parabolics.  For torsion-free
relatively hyperbolic groups with abelian parabolics, I.8(i) was
answered by the author in \cite{RelHyp} (and also follows from 
Corollary \ref{LimitHopf} above), I.8(ii) is the content of the next
theorem, and I.8(iii) is answered in Section \ref{MR} below.  We also
remark that Dahmani \cite[Theorem 0.2]{Dahmani_eq} has answered 
I.8(iv) for a class of groups
which includes torsion-free relatively hyperbolic groups with virtually
abelian parabolics.

\begin{definition}
Let $G$ be a finitely generated group.  Two systems of equations 
$\Phi$ and $\Phi'$ in finitely many variables over $G$ are {\em equivalent} if the sets of solutions of $\Phi$ and $\Phi'$ are the 
same in $G^n$ (where the number of variables is $n$). 
\end{definition}

\begin{theorem} [Compare {\cite[Theorem 1.22]{SelaHyp}}] \label{FiniteSubsystem}
Suppose that $\Gamma$ is a torsion-free relatively hyperbolic group with abelian parabolics. Then every system of equations in finitely
many variables over $\Gamma$ (without coefficients) is equivalent to a finite subsystem.
\end{theorem}
\begin{proof}
We follow the proof of \cite[Theorem 1.22]{SelaHyp}.  Let $F$ be the free group with basis the variables of $\Sigma$.  We implicitly
use the fact that a solution to $\Sigma$ over $\Gamma$
corresponds to a  homomorphism from $F$ to $\Gamma$ which
sends each of the equations in $\Sigma$ (considered as words
in $F$) to $1 \in \Gamma$.

Let $\Sigma$ be a system of equations in finitely many variables over $\Gamma$.  We iteratively construct a directed locally finite tree as follows.  Start with the first equation $\sigma_1$ in $\Sigma$, and associate with it a one relator group $G_1 = F/\langle \sigma_1 \rangle^F$.  By Proposition \ref{FiniteMax}, to $G_1$ is associated a finite number of maximal $\Gamma$--limit quotients.  Denote these
quotients of $G_1$ by $R_1, \ldots , R_m$.  Place $G_1$ at the root node of a tree, and a directed edge from $G_1$ to each $R_i$.  Note that if $G_1$ is a $\Gamma$--limit group then it is 
a maximal $\Gamma$--limit quotient of itself, and we do not need any
new vertices at this stage.

Now let $\sigma_2$ be the second equation in $\Sigma$, and consider each $R_i$ in turn.  If $\sigma_2$ represents the trivial element of $R_i$, leave it unchanged.  If $\sigma_2$ is nontrivial in $R_i$, define $\widehat{R_i} = R_i / \langle \sigma_2 \rangle^{R_i}$.  With $\widehat{R_i}$, we associate its finite collection of maximal $\Gamma$--limit quotients, and extend the locally finite tree by adding new vertices for these quotients of $R_i$, and directed edges joining $R_i$ to each of its quotients.

Continue this procedure iteratively.  By Theorem \ref{SeriesTerminates}, each branch of this locally finite tree is finite, and therefore by Konig's Lemma the entire tree is finite.  This implies that the construction of this tree terminates after finitely many steps.

Let $\Sigma'$ be the (finite) subset of $\Sigma$ consisting of those
equations $\sigma_i$ considered before the above procedure terminates.
We claim that the system $\Sigma$ is equivalent to the subsystem
$\Sigma'$.  Certainly any solution to $\Sigma$ is a solution
to $\Sigma'$.  Suppose that there is a solution to $\Sigma'$ which is
not a solution to $\Sigma$.  Then there is a homomorphism
$h \co F \to \Gamma$ which sends each element of $\Sigma'$ to $1$
but so that $h(\sigma_j) \neq 1$ for some $\sigma_j \in \Sigma 
\smallsetminus \Sigma'$.

Let $F$ be the free group with basis the variables of $\Sigma$.
Any homomorphism from $F$ to $\Gamma$ which sends $\sigma_1$ to
$1$ factors through $G_1$, and then in turn factors through one of the maximal $\Gamma$--limit quotients
$R_1, \ldots , R_m$ of $G_1$.  Similarly, any homomorphism from $F$ to $\Gamma$ which
sends $\sigma_1$ and $\sigma_2$ to $1$ either factors through one
of the $R_i$ (in case $\sigma_2$ is trivial in $R_i$) or else through one of the maximal $\Gamma$--limit
quotients of $\widehat{R_i} = R_i / \langle \sigma_2 \rangle^{R_i}$.

Arguing in this manner, we see that any homomorphism from $F$ to 
$\Gamma$ which sends all of $\Sigma'$ to $1$ factors through $G_1$
and then through some branch in the tree.  Thus we find some
terminal vertex group $V$ of this tree so that $\sigma_j \neq 1$ in $V$
(since $h(\sigma_j) \ne 1$).  This implies that the
procedure does not terminate when we consider $V$.  This
contradiction implies that all solutions to $\Sigma'$ are solutions
to $\Sigma$, as required.
\end{proof}

Guba \cite{Guba} proved the analogous theorem for free groups, whilst Sela \cite[Theorem 1.22]{SelaHyp} proved it for torsion-free hyperbolic groups.

\section{Makanin--Razborov diagrams} \label{MR}

In this final section, we describe the construction of {\em Makanin--Razborov} diagrams for $\Gamma$, which give a description of the set $\text{Hom}(G,\Gamma)$, where $G$ is an arbitrary finitely generated group.  This is analogous to the constructions in \cite[Section 5]{Sela1} and \cite[Section 1]{SelaHyp}.

Let $R$ be a freely indecomposable $\Gamma$--limit group, and let $r_1, \ldots , r_m \in R$ be a fixed generating set for $R$.  We assume that we always use the generating set $\{ r_1 , \ldots , r_m \}$ to define the length of homomorphisms, and hence to find short homomorphisms.

Following \cite[page 63]{Sela1} and \cite[page 13]{SelaHyp} we say that two proper shortening quotients $S_1, S_2$ of $R$ are {\em equivalent} if there is an isomorphism $\tau \co S_1 \to S_2$ so that the canonical quotient maps $\eta_i \co R \to S_i$, for $i = 1,2$ satisfy $\eta_2 = \tau \circ \eta_1$.  This defines an equivalence relation on the set of shortening quotients of $R$, paired with the canonical quotient maps: $\{ (S_i ,\eta_i \co R \to S_i ) \}$.

Let $SQ(R,r_1,\ldots , r_m)$ be the set of (proper) shortening quotients of $R$.  On the set $SQ(R,r_1,\ldots , r_m)$ we define a partial order as follows:  given two proper shortening quotients $S_1, S_2$ of $R$, along with canonical quotients $\eta_i \co R \to S_i$, for $i=1,2$, we say that $S_1 > S_2$ if there exists a proper epimorphism $\nu \co S_1 \to S_2$ so that $\eta_2 = \nu \circ \eta_1$.

\begin{lemma} [Compare {\cite[Lemma 1.23]{SelaHyp}}] \label{MaxShortenExist}
Let $L$ be a freely-indecomposable $\Gamma$--limit group.  Let $S_1 < S_2 < S_3 < \cdots$ (where $S_j \in SQ(L,r_1,\ldots , r_m)$) be a properly increasing sequence of (proper) shortening quotients of $L$.  Then there exists a proper shortening quotient $R \in SQ(L,r_1,\ldots , r_m)$ so that for each $j$ we have $R > S_j$.
\end{lemma}
\begin{proof}
Restricting to short homomorphisms throughout, the proof is identical to that of Proposition \ref{MaxShorten} above.  We remark that the sequence of short homomorphisms $\{ h_m \}$ obtained to find  the shortening quotient $R$ as in the proof of Proposition \ref{MaxShorten} is not contained in finitely many conjugacy classes.  This is because the sequence of kernels $\text{ker}(\eta_i)$ (where $\eta_i \co L \to S_i$ is the canonical epimorphism) is strictly decreasing.

Thus, 
$R$ is indeed a proper shortening quotient of $L$, by Proposition
\ref{ShortQProps}.
\end{proof}

\begin{lemma} [Compare {\cite[Lemma 1.24]{SelaHyp}}] \label{FinMaxShortenQ}
Let $L$ be a freely-indecomposable $\Gamma$--limit group.  The set, $SQ(L,r_1,\ldots , r_m)$, of (proper) shortening quotients of $L$ contains only finitely many equivalence classes of maximal elements with respect to the partial order.
\end{lemma}
\begin{proof}
Once again, restricting throughout to short homomorphisms, the proof is almost identical to that of Proposition \ref{FiniteMax} above.

Once again, we need to ensure that the shortening quotient $M$ obtained is a proper quotient of $L$.  This will follow from 
Proposition \ref{ShortQProps} if we can ensure that the (short) 
homomorphisms $\{ h_i \co F_d \to \Gamma \}$ obtained do not
belong to finitely many conjugacy classes.  In order to ensure this, note that for each $i \neq j$, there exists $w_{ij} \in F_d$ so that $\nu_i(w_{ij}) = 1$ but $\nu_j(w_{ij}) \neq 1$, since $R_i$ and $R_j$ are maximal and inequivalent. 

  Now 
choose  a short homomorphism $h_i \co F_d \to \Gamma$ so that (i) for words $w$ of length at most $i$ in $F_d$ we have $h_i(w) = 1$ if and only if $\nu_i(w) = 1$; (ii) $h_i$ factors through $\pi \co F_d \to G$; (iii) $h_i(w_{ij}) = 1$ for $j = 1, \ldots , i-1$; and (iv) $h_i(w_{ki}) \neq 1$ for $k = i+1, \ldots , 2i$.

Conditions (iii) and (iv) ensure that for each $i$ the set $\{ h_i, \ldots , h_{2i} \}$ have distinct kernels and so belong to different conjugacy classes.  Therefore, the $h_i$ do not belong to finitely many conjugacy
classes.

The rest of the proof is identical to that of Proposition \ref{FiniteMax} above.
\end{proof}

We can now use shortening quotients to `encode and simplify' all homomorphisms from a freely-indecomposable $\Gamma$--limit group into $\Gamma$.

\begin{proposition} [Compare {\cite[Proposition 1.25]{SelaHyp}}] \label{HomoFactor}
Suppose that $R$ is a freely-indecomposable $\Gamma$--limit group.  Let $r_1, \ldots , r_m \in R$ be a generating set for $R$, and let $M_1, \ldots , M_k$ be a set of representatives of the (finite) set of equivalence classes of maximal (proper) shortening quotients in $SQ(R,r_1,\ldots , r_m)$, equipped with the canonical quotient maps $\eta_i \co R \to M_i$, for $i = 1, \ldots , k$.

Let $h \co R \to \Gamma$ be a homomorphism which is not equivalent to an embedding.  Then there exist a (not necessarily unique) index $1 \le i \le k$, a homomorphism $\widehat{h} \co R \to \Gamma$ so that
$\widehat{h} \sim h$, and a homomorphism $h_{M_i} \co M_i \to \Gamma$ so that $\widehat{h} = h_{M_i} \circ \eta_i$.
\end{proposition}
\begin{proof}
Choose $\widehat{h} \sim h$ so that $\widehat{h}$ is short.  By
hypothesis, $\widehat{h}$ is not injective.  The constant sequence 
$\widehat{h}, \widehat{h}, \ldots$ converges into a proper shortening 
quotient $S$ of $R$.  Now, $S \cong \widehat{h}(R)$, and the 
canonical quotient map is just $\widehat{h}$.  By Lemma
\ref{FinMaxShortenQ}, there exists some $M_i$ so that $M_i > S$
or $M_i$ is equivalent to $S$. The result now follows.
\end{proof}

We use Propositions \ref{FiniteMax} and \ref{HomoFactor} to prove Theorem \ref{FactorSet}.

\begin{proof}[Proof of Theorem \ref{FactorSet}]
Let $\Gamma$ be a torsion-free relatively hyperbolic group with abelian parabolics, $G$ a freely indecomposable finitely generated group and $h \co G \to \Gamma$ a homomorphism.

The image $h(G)$ is a $\Gamma$--limit group.  By Proposition \ref{FiniteMax}, there are finitely many equivalence classes of maximal $\Gamma$--limit groups which are quotients of $G$, and $h$ factors through one of these maximal $\Gamma$--limit groups.  If $G$ is not a $\Gamma$--limit group, then these maximal $\Gamma$--limit quotients of $G$ are proper quotients, and the theorem is proved in this case.

Thus we may assume that $G$ is a $\Gamma$--limit group, and the result now follows from Proposition \ref{HomoFactor}.
\end{proof}

\begin{remark} \label{HomFreeProd}
It is reasonable in the hypothesis to restrict to freely indecomposable groups $G$, because if $G = A \ast B$ then $\text{Hom}(G,\Gamma) = \text{Hom}(A,\Gamma) \times \text{Hom}(B,\Gamma)$.  Thus a 
homomorphism $h \co G \to \Gamma$ induces a pair of 
homomorphisms, $h_A \co A \to \Gamma$ and $h_B \co B \to 
\Gamma$, and vice versa.

Hence, it is easy to understand the homomorphisms from an arbitrary finitely generated group $G$ to $\Gamma$ in terms of the sets of homomorphisms of the free factors in the Grushko decomposition of $G$ to $\Gamma$.
\end{remark}

Finally, we now construct {\em Makanin--Razborov diagrams} over 
$\Gamma$.  Let $G$ be an arbitrary finitely generated group.  
A Makanin--Razborov
diagram is a finite (directed) tree associated to a finitely generated group 
$G$ which encodes the set $\text{Hom}(G,\Gamma)$.  

There is a root vertex, labelled by $G$.  All other vertices are labelled
by $\Gamma$--limit groups.  We start with $G$.  

Suppose first that $G$ is not a
$\Gamma$--limit group, then by Proposition \ref{FiniteMax} there
are finitely many equivalence classes of maximal $\Gamma$--limit quotients of $G$.  Let $R_1, \ldots, R_s$ be a collection of representatives of these equivalence classes.  Add a new vertex for
each $R_i$, and a directed edge joining $G$ to $R_i$, labelled by
the canonical quotient map from $G$ to $R_i$.

These are the only edges emanating from $G$ (unless $G$ is a 
$\Gamma$--limit group).  All other edges will be between a pair
of $\Gamma$--limit groups.

If $G$ {\em is} a $\Gamma$--limit group, it is analysed in exactly the same way
as all of the other vertices.  We proceed with this analysis now.

There are two kinds of edges.  Suppose that the $\Gamma$--limit group
$R$ is the label of a vertex, and the Grushko decomposition of $R$
is $R = H_1 \ast \ldots \ast H_k \ast F$, where each $H_i$ is freely
indecomposable, noncyclic and finitely generated, and $F$ is a finitely generated
free group.  Then we add new vertices, one for each $H_i$ and one
for $F$, and directed edges from $R$ to each $H_i$ and to $F$. We do
not label this first kind of edge.

We now describe the second kind of edge.  Suppose that $R$ is a noncyclic freely
indecomposable $\Gamma$--limit group labelling some vertex.  By 
Lemma \ref{FinMaxShortenQ}, $R$ admits only finitely many equivalence classes of maximal proper shortening quotients.  Add
a new vertex for each equivalence class, and a directed edge from $R$
to each of these quotients, labelled by the canonical quotient map.

The vertices labelled by finitely generated free groups are terminal
vertices of the tree.

We thus build the Makanin--Razborov diagram, successively taking
Grushko decompositions and maximal proper shortening quotients,
if possible.  

This tree obviously has finite width, and by 
Theorem \ref{SeriesTerminates} each branch has finite length.  
Therefore, by Konig's Lemma, the tree is finite.

We now describe how the Makanin--Razborov diagram for $G$ encodes
the set of all homomorphisms from $G$ to $\Gamma$.  Denote the
Makanin--Razborov diagram for $G$ by $\mbox{MR}_G$.

Any homomorphism from $G$ to $\Gamma$ factors
through one of the $R_i$.  Similarly, if $L$ is the label of a vertex,
and $L$ is noncyclic and freely indecomposable, with $M_1, \ldots , M_k$ the 
maximal proper shortening quotients of $L$, then any homomorphism
from $L$ to $\Gamma$ which is not equivalent to an embedding
factors through one of the $M_i$.

For any homomorphism $h \co G \to \Gamma$ we
can find a sub-tree $\Delta_h$ of $\mbox{MR}_G$. Each edge
of $\Delta_h$ is labelled as in $\mbox{MR}_G$ (or not labelled,
as the case may be).  Each vertex $v$ has two labels:  the group $H_v$ which labels it
in $\mbox{MR}_G$, and also a homomorphism $h_v \co H_v \to \Gamma$.

The diagram $\Delta_h$ encodes a factorisation of $h$, and is constructed as follows:
The root vertex of $\Delta_h$ is the root vertex of $\mbox{MR}_G$, labelled by $G$
and $h$.

If $G$ is not a $\Gamma$--limit group, then there is some maximal $\Gamma$--limit quotient
$R_i$ of $G$ (a proper quotient), equipped with the canonical quotient $\eta \co G \to R_i$,
so that $h = h' \circ \eta$ for some homomorphism $h' \co R_i \to \Gamma$.  Add to $\Delta_h$
the edge labelled by $\eta$, and vertex labelled $R_i$ (which we label by $R_i$ and $h'$).

From $R_i$ (respectively $G$, in case $G$ is a $\Gamma$--limit group), we add the vertices 
corresponding to the Grushko decomposition.  These vertices are labelled by the 
Grushko factors, and the homomorphism induced on these free factors by $h'$ 
(respectively $h$), as in Remark \ref{HomFreeProd}.

From a free Grushko factor $F_n$, of rank $n$, the set $\text{Hom}(F_n,\Gamma)$ 
is naturally parametrised by $\Gamma^n$.  The vertex labelled by $F_n$ is a terminal vertex of $\mbox{MR}_G$, and of $\Delta_h$ also 
(whatever the other label may be).

Suppose that $v$ is a vertex of $\Delta_h$, labelled by a noncyclic freely indecomposable $\Gamma$--limit group $L$, and an associated homomorphism
$h_L \co L \to \Gamma$.  Then either 

(i)\qua $h_L$ is equivalent to an embedding $\widehat{h_L} \co L \to \Gamma$, or

(ii)\qua
$h_L$ is equivalent to some homomorphism $\widehat{h_L}$ which factors through one of 
the maximal  proper shortening quotients of $L$.  

In the first case, the associated vertex is a terminal 
vertex of $\Delta_h$ (even though
it may not be a terminal vertex of $\mbox{MR}_G$; see Remark \ref{NotTerminal} below).
In the second case, let $M$ be the maximal 
shortening quotient, $\rho_M \co L \to M$ the canonical quotient,  and $\widehat{h_L} = 
h_M \circ \rho_M$.  To $\Delta_h$ we add the edge labelled by $\rho_M$,
and the vertex from $\mbox{MR}_G$ labelled by $M$, which we label by $M$ and $h_M$.

We then take the Grushko decomposition of $M$, and proceed by repeating the above analysis.  
In this way, we obtain a tree $\Delta_h$ which encodes a factorisation of $h$.  At each vertex $v$, with labels $(G_v, h_v \co G_v \to \Gamma)$, one of the following is performed: 
\begin{enumerate}
\item Take the Grushko
decomoposition of $G_v$, and consider the
homomorphisms induced by $h_v$ on the free factors in turn.
\item Replace $h_v$ by an equivalent
homomorphism $\widehat{h_v} \co G_v \to \Gamma$ which is an injection, and stop the analysis of the branch at this vertex.
\item Replace $h_v$ by an equivalent
homomorphism $\widehat{h_v} \co G_v \to \Gamma$, where
$\widehat{h_v} = h_M \circ \eta$, for some maximal proper
shortening quotient $\eta \co G_v \to M$ of $G_v$, and proceed
by analysing $(M,h_M)$.
\item Stop the analysis of a branch at a free group, since 
$\text{Hom}(F,\Gamma)$ is easily understood.
\end{enumerate}

Note that in general the diagram $\Delta_h$ need not be unique.

\begin{remark} \label{NotTerminal}
It is worth remarking that the terminal vertices of $\Delta_h$ need 
not be terminal vertices of $\mbox{MR}_G$.  This is because
it is possible that some freely indecomposable finitely generated subgroup $H$
of $\Gamma$ is also a strict $\Gamma$--limit group which admits
a proper shortening quotient.  Thus, some homomorphisms from $H$
to $\Gamma$ are not equivalent to an injection, but some are.

This phenomena arises in the case of Makanin--Razborov diagrams
over torsion-free hyperbolic groups (as in \cite{SelaHyp}), but not
in Makanin--Razborov diagrams over free groups (because
finitely generated subgroups of free groups are free, and free
groups are always terminal vertices). 
\end{remark}

Remark \ref{NotTerminal} captures one of the key differences between 
the Makanin--Razborov diagrams for free groups (as constructed in \cite[Section 5]{Sela1})
and those constructed
in this paper for torsion-free relatively hyperbolic groups with abelian
parabolics.

There is one other difference which is worthy of remark. In the case
of free groups (and of torsion-free hyperbolic groups), every 
non-injective homomorphism from a $\Gamma$--limit group $L$ to the
target group $\Gamma$ is equivalent to one which factors through
a maximal shortening quotient of $L$.  This equivalent homomorphism
can be realised by pre-composing with an element of $\text{Mod}(L)$,
and post-composing with an inner automorphism of $\Gamma$.  Since
this conjugation does not change the kernel, we 
see that for every $h \co L \to \Gamma$ there is a maximal shortening
quotient $M$ of $L$, and $\phi \in \text{Mod}(L)$ so that $h \circ \phi$
factors through $M$.  Applying this analysis repeatedly leads
to a factorisation of any homomorphism as a composition of
modular automorphisms and canonical quotient maps, followed
finally by a product of embeddings and maps from finitely
generated free groups to the target group $\Gamma$.

If $\Gamma$ is a torsion-free relatively hyperbolic group with 
abelian parabolics,
the situation is slightly more complicated.  It is still the case that
if $L$ is a $\Gamma$--limit group then for any $h \co L \to \Gamma$
which is not equivalent to an embedding, there is a maximal shortening quotient $M$ of $L$ and $\hat{h} \sim h$ so that $\hat{h}$ factors
through $M$.  However, the equivalence between $h$ and $\hat{h}$
now involves bending moves, as well as pre-composition with an
element of $\text{Mod}(L)$ and conjugation in $\Gamma$.  Since
the bending moves can change the kernel of a map,\footnote{I thank
the referee for pointing out this important fact.} it is not necessarily
the case
that there is $\phi \in \text{Mod}(L)$ so that $h \circ \phi$ factors
through some maximal shortening quotient. Therefore, we do not 
get a factorisation analogous to that for free and torsion-free hyperbolic
groups.

The construction of Makanin--Razborov diagrams for torsion-free relatively hyperbolic groups with abelian parabolics (essentially) answers a question asked by Sela \cite[Problem I.8(iii)]{SelaProblems}. See the discussion above Theorem \ref{FiniteSubsystem}.

\end{document}